\documentclass[11pt, a4paper]{article}

\usepackage{amsmath}
\usepackage{booktabs}
\usepackage{graphicx}
\usepackage{natbib}
\usepackage{color}
\usepackage[colorlinks, citecolor=blue, urlcolor=blue]{hyperref}
\usepackage[left = 2.5cm, right = 2.5cm, bottom = 3cm, top = 2.5cm]{geometry}
\usepackage{pdfpages}
\usepackage{multirow}
\usepackage{authblk}

\newcommand{\bRe}{{\overline{Re}}}
\newcommand{\sumq}{{\sum_{i=1}^q}}

\newcommand{\bBC}{{\overline{BC}}}
\def\P{{\mathop{\rm pr}}}
\def\E{{\mathop{\rm E}}}
\def\Cov{{\mathop{\rm Cov}}}
\def\V{{\mathop{\rm Var}}}

\author[1]{Ruggero Bellio\thanks{ruggero.bellio@uniud.it}}
\author[2, 3]{Ioannis Kosmidis\thanks{ioannis.kosmidis@warwick.ac.uk}}
\author[4]{Alessandra Salvan\thanks{alessandra.salvan@unipd.it}}
\author[4]{Nicola Sartori\thanks{nicola.sartori@unipd.it}}

\affil[1]{Department of Economics and Statistics, University of Udine, Udine, 33100, Italy}
\affil[2]{Department of Statistics, University of Warwick, Coventry, CV4 7AL, UK}
\affil[3]{The Alan Turing Institute, London, NW1 2DB, UK}
\affil[4]{Department of Statistical Sciences, University of Padova, 35121, Padova, Italy}

\title{Parametric bootstrap inference for stratified models with high-dimensional nuisance specifications}

\begin{document}
\maketitle

\begin{abstract}
Inference about a scalar parameter of interest typically relies on the
asymptotic normality of common likelihood pivots, such as the signed
likelihood root, the score and Wald statistics. Nevertheless, the
resulting inferential procedures are known to perform poorly when the
dimension of the nuisance parameter is large relative to the sample
size and when the information about the parameters is limited. In many
such cases, the use of asymptotic normality of analytical
modifications of the signed likelihood root is known to recover
inferential performance. It is proved here that parametric bootstrap
of standard likelihood pivots results in as accurate inferences as
analytical modifications of the signed likelihood root do in
stratified models with stratum specific nuisance parameters. We focus
on the challenging case where the number of strata increases as fast
or faster than the stratum samples size.  It is also shown that this
equivalence holds regardless of whether constrained or unconstrained
bootstrap is used. This is in contrast to when the number of strata is
fixed or increases slower than the stratum sample size, where we show that
constrained bootstrap corrects inference to a higher order than
unconstrained bootstrap. Simulation experiments support the
theoretical findings and demonstrate the excellent performance of
bootstrap in extreme scenarios.
\end{abstract}

\section{Introduction}
Standard likelihood inference about a scalar parameter of interest is
based on the asymptotic normality of likelihood pivots, such as the
signed likelihood root, score and Wald statistics. This asymptotic
approximation can be quite inaccurate in the presence of many nuisance
parameters. An alternative, which guarantees higher accuracy, is based
on the asymptotic normality of analytical modifications of the
signed likelihood root, generally termed modified signed likelihood root
\citep[see for instance][Chapter 7]{Severini:2000}. In a two-index
stratified asymptotic setting, in which both the dimension of the data
and the number of nuisance parameters grow, the modified signed
likelihood root has been proved to be highly accurate even in rather
extreme scenarios with many nuisance parameters and very limited information
\citep{Sartori:2003}.

Parametric bootstrap methods provide an alternative assessment of tail
probabilities for likelihood pivots and, in standard asymptotic
settings, where the number of nuisance parameters is fixed and
regularity conditions are satisfied \citep[Section
3.4]{Severini:2000}, have been shown to guarantee an equivalent level of
asymptotic accuracy as analytical modifications of the signed
likelihood root \citep[see][Chapter 11]{Young.Smith:2005}. In
particular, the two main variants of parametric bootstrap are
constrained and unconstrained bootstrap (also know as conventional
bootstrap). In the latter, the sampling distribution of the statistic
is computed at the full maximum likelihood estimate, and in the former
at the constrained maximum likelihood estimate for a given value of
the parameter of interest. In standard asymptotic settings,
constrained bootstrap \citep{DiCiccio.etal:2001, Lee.Young:2005}
corrects inference about a scalar parameter in the presence of
nuisance parameters to a higher order than unconstrained bootstrap. On the other hand, numerical differences are rarely detectable.
Although bootstrap methods are, typically, more computationally
demanding than analytical approximations to the distribution of
pivots, they are available in some non-regular cases in which the
modified signed likelihood root is not computable.

We investigate the properties of parametric bootstrap in models for
stratified data in a two-index asymptotic setting, where both the
number $q$ of strata and the sample size $m$ of each stratum grow. In
this setting, the usual likelihood pivots are asymptotically standard
normal provided $q=o(m)$, while the condition for the modified
signed likelihood root
 is $q=o(m^3)$ \citep{Sartori:2003}. If
$q = O(m^\alpha)$, then for $0 \leq \alpha < 1$ the asymptotic
normality of standard likelihood pivots still holds, with error of
order $O_p(m^{(\alpha-1)/2})$ \citep[][formula (8)]{Sartori:2003},
while asymptotic normality fails in the highly stratified case with
$\alpha \geq 1$.
In that case, the aim of higher-order solutions is to recover
first-order validity of inferential procedures.

We show here that parametric bootstrap provides valid inference when
$q = O(m^\alpha)$, provided that $\alpha < 3$. In particular, if
$0 \leq \alpha < 1$, constrained bootstrap is theoretically more
accurate than unconstrained bootstrap, and both improve over standard
first-order asymptotic results. On the other hand, when
$1 \leq \alpha < 3$
both variants of parametric bootstrap are equally accurate, recovering
first-order accuracy with the same order of error as higher-order
analytical solutions.

The theoretical results are
supported by extensive simulation studies, which illustrate that
parametric bootstrap is at least as accurate as use of the modified
signed likelihood root, and provide evidence that constrained bootstrap  can be even more
accurate in some very extreme scenarios.

\section{Background}
\label{sec:inference}
Let $l(\theta; y)$ be the log-likelihood function for a parameter
$\theta$ based on a sample $y$ of size $n$, which is considered to be a
realization of a random vector $Y$. We  treat the case where the
vector of parameters is partitioned as
$\theta=(\psi,\lambda^\top)^\top$, where $\psi$ is a scalar parameter
of interest and $\lambda$ is a vector of nuisance parameters, and
denote by $\hat\theta(y)=(\hat\psi(y), \hat\lambda(y)^\top)^\top$ the
maximum likelihood estimate of $\theta$ and by
$\hat\theta_\psi(y)=(\psi,\hat\lambda_{\psi}(y)^\top)^\top$ the
constrained maximum likelihood estimate of $\theta$ for fixed
$\psi$. We let $U(\theta; y)=\nabla l(\theta; y)$ denote the score
vector, and 
$j(\theta; y)=-\nabla\nabla^\top l(\theta; y)$ the observed
information, with $i(\theta)= \E_\theta\{j(\theta; Y)\}$ denoting the expected
information. The argument $\theta$ will be dropped when no ambiguity
arises, and components of vectors and blocks of matrices will be
denoted by subscripts, so that for instance $U_\psi(\theta;y)$ denotes the component of the score vector corresponding to $\psi$. 
Furthermore, the argument $y$ will be dropped
whenever evaluation is at the random variable $Y$ instead of the
sample $y$. For example, $U_\psi = U_\psi(\theta; Y)$,
$U_\lambda = U_\lambda(\theta; Y)$, $i_{\psi\psi} = i_{\psi\psi}(\theta)$
and $i_{\psi\lambda}=i_{\psi\lambda}(\theta)$ are the $(\psi, \psi)$
and $(\psi, \lambda)$ blocks of $i(\theta)$, and so on.

The signed likelihood root, the score statistic and Wald statistic for
inference about $\psi$ are
\begin{eqnarray}
  R(\psi; y) &=& \mathop{\rm sign}\left(\hat\psi(y) -
                 \psi\right)\sqrt{2\left\{l(\hat\theta(y); y) -
                 l(\hat\theta_\psi(y); y)\right\} }\,\label{R} \\
  S(\psi; y) &=& \frac{U_p(\psi; y)}{\sqrt{i_{\psi\psi |\lambda}(\hat\theta_\psi(y))}}\,,\label{score}\\
  T(\psi; y) &=&   (\hat\psi(y) - \psi)\sqrt{j_p(\hat\psi(y); y)} \, , \label{wald}
\end{eqnarray}
respectively, where $U_p(\psi; y)=U_\psi(\hat\theta_\psi(y); y)$ is
the profile score, $j_p(\psi; y)=- d U_p(\psi; y) / d \psi$ is the profile
observed information and
$i_{\psi\psi
  |\lambda}=i_{\psi\psi}-i_{\psi\lambda}i_{\lambda\lambda}^{-1}i_{\lambda\psi}$
is the partial information about $\psi$. While (\ref{R}) and
(\ref{score}) are invariant with respect to interest respecting
reparameterizations, (\ref{wald}) is not.

Computation of $p$-values and confidence intervals for $\psi$ requires
the distribution of statistics (\ref{R}), (\ref{score}) and
(\ref{wald}). In standard asymptotic settings, one possibility is to
rely on the first-order asymptotic normal approximation to their
distribution.  For instance,
$\P_\theta\{R(\psi) \leq R(\psi; y)\}= \Phi
(R(\psi; y))\{1+O(n^{-1/2})\}$, where $\Phi(\cdot)$ denotes the
standard normal distribution function. Improved accuracy can
be obtained with higher-order modifications $R^*(\psi; y)$ of $R(\psi; y)$,
such that
$\P_\theta\{R(\psi) \leq R(\psi; y)\}= \Phi
(R^*(\psi; y))\{1+O(n^{-1})\}$. \cite{Barndorff-Nielsen:1986}
developed a modified signed likelihood root $R^*(\psi)$ which is
standard normal with error of order $O(n^{-3/2})$.  Following this
seminal work, there have been various alternative versions of
$R^*(\psi; y)$ \citep[see][for an accessible
overview]{Pierce.Bellio:2017}.

An alternative to the asymptotic approximations to the distribution of
(\ref{R}), (\ref{score}) and (\ref{wald}) is parametric bootstrap,
which provides higher-order approximations for $p$-values, such as
$\P_\theta\{R(\psi) \leq R(\psi; y)\}$. There are two main variants of
parametric bootstrap: i) unconstrained bootstrap where samples are
simulated from the model at $\hat\theta(y)$, and ii) constrained
bootstrap where samples are simulated at $\hat\theta_\psi(y)$
\citep[see][Chapter
11]{DiCiccio.etal:2001,Lee.Young:2005,Young.Smith:2005}.

In standard asymptotic settings, unconstrained bootstrap provides
second-order accuracy. Let 
$G_{\theta}(\cdot)$ denote the distribution function of $R(\psi)$ at $\theta$, so that $G_{\theta}(R(\psi)) $ is exactly
uniform. If data $y^k$ are simulated from the model with
parameter $\hat\theta(y)$, $k=1,\ldots,K$, 
then  $p$-values for
(\ref{R}) calculated as
\begin{equation}\label{boot1}
\hat{p}^{R}_1(\psi) = \frac{1}{K} \sum_{k=1}^{K} I\{ R(\hat\psi(y); y^k) \leq R(\psi; y)\}
\end{equation}
are Monte Carlo estimates of
 $G_{\hat\theta}(R(\psi))$, which is  uniform on $(0,1)$ under repeated sampling with error of order $O(n^{-1})$, i.e.
\begin{equation}
  \label{inttrasf2}
  \P_{\theta}\left(G_{\hat\theta}(R(\psi)) \leq u \right) = u+O(n^{-1})\,.
\end{equation} 
In (\ref{boot1}), $I\{\cdot\}$ is the indicator function.

In contrast,
constrained bootstrap provides third-order accuracy;  if data $y^k$ are
simulated at $\hat\theta_\psi(y)$,  $k=1,\ldots,K$, $p$-values for
(\ref{R}) calculated as
\begin{equation}\label{boot2}
\hat{p}^{R}_2(\psi) = \frac{1}{K} \sum_{k=1}^{K} I\{ R(\psi; y^k) \leq R(\psi; y)\}
\end{equation}
are Monte Carlo estimates of $G_{\hat\theta_\psi}(R(\psi))$, which is
uniform on $(0,1)$ under repeated sampling with error of order
$O(n^{-3/2})$ \citep{Lee.Young:2005}, i.e.
\begin{equation} \label{inttrasf1}
  \P_{\theta}\left(G_{\hat\theta_\psi}(R(\psi)) \leq u \right) =u+O(n^{-3/2})\,.
\end{equation}
Similar results hold for $S(\psi)$ and $T(\psi)$
\citep{Lee.Young:2005,Young:2009} with $p$-values $\hat{p}^{S}_1$ and
$\hat{p}^{S}_2$, and $\hat{p}^{T}_1$ and $\hat{p}^{T}_2$, respectively.

As \citet[][Section 11.4]{Young.Smith:2005} note, the theoretical
advantage of constrained over unconstrained bootstrap is rarely
supported by numerical evidence, because both types of bootstrap are
able to equally improve over first-order results.

The advantage of bootstrap $p$-values in (\ref{boot1}) and
(\ref{boot2}) over the use of analytical modifications to common
statistics is that bootstrap does not require any additional, often
tedious, algebraic derivation and implementation of the necessary
modifications. Moreover, there are non-standard modelling settings,
where $R(\psi; y)$ is computable while $R^*(\psi; y)$ is not. One instance
is when one or more components of $\hat\theta(y)$ are on the
boundary of the parameter space.  The main disadvantage of bootstrap
is the additional computation that is typically required for the
repeated model fits, which can be partly mitigated by parallel
computing.

In some special cases, the distribution of (\ref{R}), (\ref{score})
and (\ref{wald}) depends only on $\psi$, so that constrained
bootstrap, as well as simulating data at
$(\psi, \hat\lambda(y)^\top)^\top$ or even at $(\psi, \lambda_0^\top)^\top$ for
arbitrary nuisance vectors $\lambda_0$, produces samples from the hypothesized
model.  This is the case when the model for fixed $\psi$ is a
transformation model \citep[see][Section 1.3]{Severini:2000}. For
instance, if $y$ is a realization of $Y=(Y_1,\ldots,Y_n)^\top$ with
independent and identically distributed components from a shape and scale model with generic
density
\[
g(y_{i};\psi, \lambda)=\frac{1}{\lambda}g^0(y_{i}/\lambda;\psi)\,,
\]
we may write $Y_{i}=\lambda Y_{i}^0$, with
$Y_{i}^0\sim g^0(y_{i};\psi)=g(y_{i};\psi,1)$. Hence, due to
equivariance of the maximum likelihood estimator, $\hat\lambda$ and
$\lambda \hat\lambda^0$ have the same distribution, where
$\hat\lambda^0$ is the maximum likelihood estimator of $\lambda$ based
on $Y_{i}^0$'s. The same representation holds for
$\hat\lambda_{\psi}$, so that the profile likelihood ratio
\[
\exp\{l(\hat\psi,\hat\lambda)- l(\psi,\hat\lambda_{\psi})\}=\prod_{i=1}^n \frac{\hat\lambda_{\psi}g^0(Y_{i}/\hat\lambda;\hat\psi)}{\hat\lambda g^0(Y_{j}/\hat\lambda_{\psi};\psi)}
\]
has the same distribution as
\[
\prod_{i=1}^n 
\frac
{\lambda\hat\lambda_{\psi}^0 g^0\left(\frac{\lambda Y_{i}^0}{\lambda\hat\lambda^0};\hat\psi\right)}
{\lambda\hat\lambda^0 g^0\left(\frac{\lambda Y_{i}^0}{\lambda\hat\lambda_{\psi}};\psi\right)}=\prod_{i=1}^n 
\frac
{ \hat\lambda_{\psi}^0 g^0(Y_{i}^0/ \hat\lambda^0;\hat\psi)}
{ \hat\lambda^0 g^0(Y_{i}^0 /\hat\lambda_{\psi};\psi)}\,,
\]
which depends on $\psi$ only.

An example with a stratified gamma model is provided in the
Supplementary Materials where simulation results confirm the exactness
of the constrained bootstrap.

\section{Two-index asymptotic theory for stratified models}

We consider a stratified setting with $q$ independent strata with $m$
observations each.  Therefore, the total number of observations is
$n=mq$.  The models considered here have
$\lambda=(\lambda_1,\ldots,\lambda_q)^\top$ as nuisance parameter,
where $\lambda_i$ is a stratum-specific parameter.  Let
$y_i=(y_{i1},\ldots, y_{im})^\top$, $i=1,\ldots,q$, denote the vector
of observations in the $i$th stratum and let
$y=(y_1^\top,\ldots,y_q^\top)^\top$.  The vectors $y_1, \ldots, y_q$
are assumed to be realizations of independent random variables
$Y_{1}, \ldots, Y_{q}$ from a parametric model with densities
$g_1(y_{1}; \psi, \lambda_1),\ldots, g_q(y_{q}; \psi, \lambda_q)$,
respectively.
The observations within strata are also assumed to be realizations of
independent random variables, so that
$g_i(y_{i}; \psi, \lambda_i) =\prod_{j=1}^m g_{ij}(y_{ij}; \psi,
\lambda_i) $, where $g_{ij}(\cdot)$ may be conditional on a covariate
vector $x_{ij}$.  Under this specification, for fixed $\psi$, the
likelihood has separable parameters $\lambda_1, \ldots, \lambda_q$, so
that $U_p(\psi) = \sumq U^i_\psi(\psi, \hat\lambda_{i\psi})$, where
$U^i_\psi$ is the contribution to $U_\psi$ from the  $i$th stratum.

We work in a two-index asymptotic setting where $q$ increases with
$m$, as $q=O(m^\alpha)$, $\alpha > 0$. The case $\alpha = 0$ corresponds
to the standard asymptotic setting. \citet[][Section 4]{Sartori:2003}
showed that $R(\psi)$, $S(\psi)$
and $T(\psi)$ are asymptotically equivalent to order $o_p(1)$ for $\alpha\geq 0$. 
Specifically, when $0\leq\alpha<1$, the equivalence of the three quantities holds with relative error of order $O_p(n^{-1/2})=O_p(m^{-(\alpha+1)/2})$, and these are asymptotically standard normal. On the other hand,  when $\alpha \geq 1$, asymptotic equivalence of $R(\psi)$, $S(\psi)$ and $T(\psi)$ holds with error of order  $O_p(m^{-1})$ and, more critically, the three statistics are not asymptotically
standard normal, so that, for instance, $\Phi\{R(\psi)\}$ is not asymptotically
uniform.

The derivation of the results is more straightforward for $S(\psi)$ because the profile score is the sum of strata profile
scores. However, the same results hold also for $R(\psi)$ and $T(\psi)$, since, as recalled above, they are both asymptotically equivalent to $S(\psi)$.  
Let $F_\theta(\cdot)$ denote the distribution function of $S(\psi)$ under $\theta$, so that $F_{\theta}(S(\psi)) $ is exactly
uniform.

The core result of the paper is that asymptotic validity of
both constrained and unconstrained bootstrap is guaranteed even in a
two-index asymptotic setting, provided that
$\alpha < 3$, that is $q = o(m^3)$. The latter condition is the same
as the one required for validity of inference based on the modified
signed likelihood root $R^*(\psi)$ \citep{Sartori:2003}. In
particular, we show that, when $0 < \alpha < 1$,
\begin{equation}
  \label{res1.l1}
\P_{\theta}\left(F_{\hat\theta_\psi}(S(\psi)) \leq u \right)
= u+O(m^{(\alpha-3)/2})
\end{equation}
and 
\begin{equation}
  \label{res2.l1}
\P_{\theta}\left(F_{\hat\theta}(S(\psi)) \leq u \right)
= u+O(m^{-1})\,,
\end{equation} 
while, when $1 \leq \alpha < 3$,
\begin{equation}
  \label{res1}
\P_{\theta}\left(F_{\hat\theta_\psi}(S(\psi)) \leq u \right)
= u+O(m^{(\alpha-3)/2})
\end{equation}
and 
\begin{equation}
  \label{res2}
\P_{\theta}\left(F_{\hat\theta}(S(\psi)) \leq u \right)
= u+O(m^{(\alpha-3)/2})\,.
\end{equation} 
Hence, when $1 \leq \alpha < 3$, the same order of error is obtained both with constrained and
unconstrained bootstrap, unlike what happens with $0\leq \alpha<1$.
The case $\alpha=0$ corresponds to the standard asymptotic setting  in which $n=O(m)$, and (\ref{res1.l1}) and (\ref{res2.l1}) reduce to
(\ref{inttrasf1}) and (\ref{inttrasf2}), respectively. A first intuition about why the two types of bootstrap have the same accuracy when $\alpha \geq 1$ is that the major effect of both bootstrap procedures is to remove the diverging bias term of the statistic, which overshadows any minor differences in theoretical performance that are found when $0\leq \alpha <1$. A formal  development of the result  is given below.

In the following we will concentrate on the more extreme case, i.e. $\alpha \geq 1$, while the proof of (\ref{res1.l1}) and (\ref{res2.l1})  for the case $0 < \alpha < 1$ is given in the Supplementary Materials.
In order to prove both (\ref{res1}) and (\ref{res2}) we need some
preliminary results about the distribution function $F_\theta(x)$ of
$S(\psi)$ in the two-index asymptotic setting.  From \citet[][formula
(6)]{Sartori:2003}, $U_p = U_p(\psi)$ can be expanded as
\begin{equation}\label{ns2003}
U_p = U_{\psi|\lambda} + B + Re\,,
\end{equation}
where $U_{\psi|\lambda}= U_\psi - i_{\psi\lambda} i_{\lambda\lambda}^{-1} U_\lambda=O_p(\sqrt{n})=O_p(m^{(\alpha+1)/2})$, having zero mean and variance   
$i_{\psi\psi |\lambda}$,  
$B=B(\theta)=O_p(m^{\alpha})$ and, with $\alpha>1$,  $Re=O_p(m^{\alpha-1})$. Details about the orders in (\ref{ns2003}) are provided in the Appendix. When $0\leq \alpha <1$ the terms in (\ref{ns2003}) are in descending order. Instead, when $1 \leq \alpha < 3$, $B$ becomes the leading term, followed by $U_{\psi|\lambda}$. Finally, $U_{\psi|\lambda}$ is dominated by $Re$  as well when $\alpha \geq 3$. In practice, when $1 \leq \alpha < 3$, bootstrap procedures, as well as higher-order analytical solutions, are able to correct for  $B$, so that $U_{\psi|\lambda}$ is again the leading term in expansion (\ref{ns2003}).

Let $M(\theta) = \E_\theta(S(\psi))$ and $\V_\theta(S(\psi))$ be the
expectation and variance of $S(\psi)$. Asymptotic expansions detailed
in the Appendix can be used to show that
\begin{eqnarray}
M(\theta) & = & 
\frac{b(\theta)}{i_{\psi\psi |\lambda}(\theta)^{1/2}}+ M_1(\theta)+ O(m^{(\alpha-5)/2}) \label{ESexp} \\
\V_\theta(S(\psi)) & = & 1 + v(\theta) + O(1/m^2) \label{VSexp}\,,
\end{eqnarray}
where $b(\theta)= \E_\theta (B)=O(m^\alpha)$, $M_1(\theta)=O(m^{(\alpha-3)/2})$ and  $v(\theta)=(\V_\theta(B) + 2\,\E_\theta(U_{\psi|\lambda} B))/i_{\psi\psi |\lambda}=O(m^{-1})$. 
The cumulants of $S(\psi)$ of order $r \in \{3, 4, \ldots\}$ are
$O\left( m^{(\alpha+1)(1-r/2)} \right) = O(n^{1-r/2})$, as in standard
asymptotics.

For the development here, we assume that the distribution function of
$S(\psi)$ admits a valid Edgeworth expansion. \citet[][Sections
5.1-5.4]{Severini:2000} gives conditions and details for the extension
of Edgeworth expansions for independent and identically distributed
random variables to likelihood pivots, such as $R(\psi)$, $S(\psi)$,
$T(\psi)$. The basic requirement, in the continuous case, is that an
Edgeworth expansion exists for the joint distribution of
log-likelihood derivatives up to the third order, implying
\begin{equation}\label{Etrue}
F_\theta(x) = \P_{\theta}\left(S(\psi)\leq x \right) = \Phi\left(\frac{x-M(\theta)}{\sqrt{\V_\theta(S(\psi))}}\right) + O(m^{-(\alpha+1)/2}) \,,
\end{equation} 
where the order of the remainder term is that of the third cumulant of $S(\psi)$.
Let $x^*(\theta)= (x-M(\theta))/\sqrt{1+ v(\theta)}$. Then
\[
\frac{x-M(\theta)}{\sqrt{\V_\theta(S(\psi))}}=
\frac{x-M(\theta)}{\sqrt{1+v(\theta)+O(m^{-2})}}=
x^*(\theta)+O(m^{-2})
\]
and 
\begin{equation}\label{E1}
F_\theta(x)= \Phi\left(x^*(\theta)\right)+ O\left( m^{-\min\left(2, \frac{\alpha+1}{2}\right)} \right)\,.
\end{equation}

We first focus on constrained bootstrap. From (\ref{E1}), 
\begin{equation}\label{E3}
F_{\hat\theta_\psi}(x)= \Phi\left(x^*(\hat\theta_\psi)\right)+O_p\left( m^{-\min\left(2, \frac{\alpha+1}{2}\right)} \right)\,.
\end{equation}
The Taylor expansions in the Appendix give 
 \begin{equation}\label{BS}
M(\hat\theta_\psi)= M(\theta)+\Delta+O_p\left(m^{-\min\left(1, \frac{5-\alpha}{2}\right)}\right)
\end{equation}
and 
\begin{equation}\label{v}
 v(\hat\theta_\psi)=v(\theta)+ O_p(m^{-2})\,,
\end{equation}
where $\Delta =O_p(m^{(\alpha-3)/2})$ is given in expression
(\ref{Delta}) of the Appendix.  Using (\ref{BS}) and (\ref{v}), we can
write
$x^*({\hat\theta_\psi})=x^*(\theta) - \Delta + O_p(m^{-\min(1, (5-\alpha)/2})$.
As a result, if $\alpha < 3$, then the following Taylor expansion of
(\ref{E3}) holds
\begin{equation}\label{CBE}
F_{\hat\theta_\psi}(x)=
F_{\theta}(x)-\phi(x^*(\theta)) \Delta+O_p(m^{-1})\,,
\end{equation}
where the error is of order $O_p(m^{-1})$ because, for $\alpha < 3$,  $\min\left(1, (5-\alpha)/2\right)=1$, while the error term in
(\ref{E3}) is $o_p(m^{-1})$ whenever $\alpha > 1$.

In order to prove (\ref{res1}), note that
$F_{\hat\theta_\psi}(S(\psi)) \leq u$ is equivalent to
$S(\psi) \leq s_u$, with $s_u$ the $u$-quantile of
$F_{\hat\theta_\psi}(\cdot)$, such that
$F_{\hat\theta_\psi}(s_u)=u$. Let $s_u^0$ be the $u$-quantile of
$F_\theta(\cdot)$. It is useful to express $s_u$ in terms of
$s^0_u$. Using (\ref{CBE}),
\[
u=F_\theta(s_u^0)=F_{\hat\theta_\psi}(s_u)=
F_{\theta}(s_u)-\phi(s_u^*(\theta)) \Delta +O_p(m^{-1})\,,
\]
where $s_u^*(\theta)=(s_u-M(\theta))/\sqrt{1+ v(\theta)}$.  Hence,
$F_{\theta}(s_u) - F_{\theta}(s_u^0) = \phi(s_u^*(\theta)) \Delta + O_p(m^{-1})$.
On the other hand, letting $F'_\theta(x) = dF_\theta(x)/dx$, from
\[
F_{\theta}(s_u^0)=F_{\theta}(s_u)+(s_u^0-s_u)F'_\theta(s_u)+O_p((s_u^0-s_u)^2)
\]
and
\[
F'_\theta(x)=\phi(x^*(\theta))/\sqrt{1+v(\theta)}+O(m^{-(\alpha+1)/2})= \phi(x^*(\theta)) +O(m^{-1})
\]
we get
\[
s_u=s_u^0+\Delta+O_p(m^{-1})+O_p(m^{\alpha-3})\,,
\]
where the $O_p(m^{\alpha-3})$ term on the right hand side comes from
$O_p((s_u^0-s_u)^2)$.  Hence, $S(\psi) \leq s_u$ is equivalent to
$S(\psi)\leq s_u^0+\Delta+O_p(m^{-1})+O_p(m^{\alpha-3})$, and
\[
\P_{\theta}\left(F_{\hat\theta_\psi}(S(\psi)) \leq u \right)=
\P_{\theta}\left(\bar{S}(\psi) \leq F_{\theta}^{-1}(u)\right)\,,
\]
where $\bar{S}(\psi)=S(\psi) -\Delta+O_p(m^{\alpha-3})+O_p(m^{-1})$, with $\Delta$ given by (\ref{Delta}), and such that  
$\E_\theta(\Delta)= O(m^{(\alpha-3)/2})$. Moreover, we have
\begin{eqnarray}
\E_\theta(\bar{S}(\psi)) &=& \E_\theta(S(\psi)) +O(m^{(\alpha-3)/2}) \,, \label{EVSprime}\\
\V_\theta(\bar{S}(\psi))&=& \V_\theta(S(\psi)-\Delta) + O(m^{-2})\, \nonumber \\
&=& \V_\theta(S(\psi)) +\V_\theta(\Delta) - 2 \Cov_\theta(S(\psi),\Delta) + O(m^{-2}) \nonumber \\
&=& \V_\theta(S(\psi)) + O(m^{-2}) \label{VarSprime} \,,
\end{eqnarray}
since $\V_\theta(\Delta)= O(m^{-2})$ and $\Cov_\theta(S(\psi),\Delta) = O(m^{-2})$, where the order of the latter is determined by the orthogonality between $U_{\psi|\lambda}$ and the leading term of $b_1(\theta)$ in (\ref{b1}).
Finally, (\ref{res1}) holds because
\begin{eqnarray*}
\P_{\theta}\left(\bar{S}(\psi) \leq F_{\theta}^{-1}(u)\right)&=& 
\P_{\theta}\left(S(\psi) \leq F_{\theta}^{-1}(u)\right)+O(m^{(\alpha-3)/2}) + O(m^{-2}) \\
&=& \P_{\theta}\left(S(\psi) \leq F_{\theta}^{-1}(u)\right)+O(m^{(\alpha-3)/2})\,\\
&=&u+O(m^{(\alpha-3)/2})\,.
\end{eqnarray*}

The proof of (\ref{res2}) for unconstrained bootstrap is obtained
along the same steps as above. 
In particular, expansion (\ref{BE}) holds for 
$F_{\hat\theta}(x)$, having the same form as (\ref{CBE}), with $\Delta$ replaced by $\Delta_1$, which is still of order $O_p(m^{(\alpha-3)/2})$. Details are given in the Appendix.
However, while (\ref{EVSprime}) is still true, (\ref{VarSprime}) holds with an error of order $O(m^{-1})$, because there is no orthogonality between $U_{\psi|\lambda}$ and the leading terms of $b_2(\theta)$, given in (\ref{expBPconv}). Therefore, for unconstrained bootstrap we have
\begin{eqnarray*}
\P_{\theta}\left(\bar{S}(\psi) \leq F_{\theta}^{-1}(u)\right)&=& 
\P_{\theta}\left(S(\psi) \leq F_{\theta}^{-1}(u)\right)+O(m^{(\alpha-3)/2}) + O(m^{-1}) \\
&=& \P_{\theta}\left(S(\psi) \leq F_{\theta}^{-1}(u)\right)+O(m^{(\alpha-3)/2})\,\\
&=&u+O(m^{(\alpha-3)/2})\,.
\end{eqnarray*}
Hence, when $\alpha \geq 1$, errors in (\ref{res1}) and (\ref{res2})
are of the same order because the $O(m^{(\alpha-3)/2})$ error in the mean of $\bar S(\psi)$ dominates the $O(m^{-2})$ and $O(m^{-1})$ errors in the variance of $\bar S(\psi)$ in the constrained and unconstrained cases, respectively. However, the different errors in the variance of $\bar S(\psi)$ may have some effects and explain why the constrained bootstrap is sometimes numerically more accurate in extreme settings.

The arguments used in the proofs of~(\ref{res1}) and~(\ref{res2})
suggest that a location and scale adjustment to the statistic, as done
for $R(\psi)$ in a standard asymptotic setting by
\citet{DiCiccio.etal:2001} and \citet{Stern:2006}, is the key
requirement to recover approximate uniformity of $p$-values. In this
respect, a bootstrap location and scale adjustment of $R(\psi)$,
$S(\psi)$ or $T(\psi)$ is expected to be as effective as bootstrapping
the distribution of the statistic. This conjecture is confirmed by the
numerical results, both in the following section and in the
Supplementary Materials.

%%%%%%%%%%%%%%%%%%%%%%%%%%%%%%%%%%%%%%%%%%%%%%%%%%%%%%%%%%%%%%%%%%%%%%%%%%%%%%%%%%%%%%%%%%%%%%%%%%%%%%%%%%%%%%%%%%%%%%%%%%%%

\section{Simulation studies}
\label{sec:simulations}

The finite-sample properties of unconstrained and constrained
parametric bootstrap are assessed through extensive simulation
studies, for three statistical models for stratified data. In particular, we consider 
a beta model, a curved exponential family model and a truncated regression model,
with the results for further models reported in the Supplementary Materials.
 For each model,
we conduct 9 simulation experiments, one for each combination of
number of strata $q \in \{10, 100, 1000\}$ and stratum sample size
$m \in \{4, 8, 16\}$.

Each simulation experiment involves $10000$ simulated samples under
the model at a fixed parameter vector
$\theta_0 = (\psi_0, \lambda_0^\top)^\top$. For each simulated sample,
17 statistics and 6 bootstrap-based $p$-values are computed for
testing $\psi = \psi_0$. In particular, the statistics that are
computed are i) $R(\psi)$, $S(\psi)$, $T(\psi)$, ii) the location and
location-and-scale adjusted versions of $R(\psi)$, $S(\psi)$,
$T(\psi)$, where the mean and variance of each statistic are estimated
using unconstrained bootstrap (at $\hat\theta$) and constrained
boostrap (at $\hat\theta_\psi$), and iii) $R^*(\phi)$ and the signed
likelihood root computed from the modified profile likelihood
\citep[see, for instance][Chapter 8]{Severini:2000}. The higher-order
adjustment required for the latter two statistics is obtained using
expected moments of likelihood quantities as in \citet[Section
7.5]{Severini:2000}. Finally, for each of $R(\psi)$, $S(\psi)$, and
$T(\psi)$, we compute the unconstrained and constrained bootstrap
$p$-values in~(\ref{boot1}) and in~(\ref{boot2}), respectively.

In the interest of space, in what follows, we only report results for
the 6 statistics based on $R(\psi)$ shown in
Table~\ref{tab:stats}. The conclusions for the remaining statistics
and $p$-values are qualitatively the same. Results are also only
presented for $(q, m) = (10, 4)$, $(q, m) = (100, 4)$,
$(q, m) = (1000, 4)$, $(q, m) = (1000, 8)$, and $(q, m) = (1000, 16)$,
because these combinations of $q$ and $m$ are sufficient for assessing
the performance of the statistics as $q$ and $m$ grow. The results
from all simulation experiments are provided in the Supplementary
Materials.

\begin{table}[t!]
  \caption{Statistics considered for the results of the simulation
    experiments. The mean $\tilde{\mu}^{R}$ and the standard deviation
    $\tilde{\sigma}^{R}$ of $R(\psi)$ are estimated through
    constrained bootstrap, by simulating from the model at
    $\theta = \hat\theta_{\psi}$.}
  \begin{center}
      \begin{tabular}{p{0.23\textwidth}p{0.12\textwidth}p{0.45\textwidth}}
        \toprule
        Statistic & Plotting  &  Description   \\  
          &  Symbol &    \\
        \midrule
        $R(\psi)$ &  $R$ & Signed likelihood root \\ 
        $R^{*}(\psi)$ &  $R^{*}$ & Modified signed likelihood root \\ 
        $\Phi^{-1}\{\hat{p}^{R}_1(\psi)\}$ &  $R^{u}$  & Transformed $p$-value from unconstrained bootstrap  of $R(\psi)$\\ 
        $\Phi^{-1}\{\hat{p}^{R}_{2}(\psi)\}$ &  $R^{c}$ & Transformed $p$-value from constrained bootstrap  of $R(\psi)$\\ 
        $R(\psi) - \tilde{\mu}^{R}$ & $R_{l}^{c}$ & Location adjusted $R(\psi)$   \\
        $(R(\psi) - \tilde{\mu}^{R}) / \tilde{\sigma}^{R}$ & $R_{ls}^{c}$ & Location-and-scale adjusted $R(\psi)$  \\ \bottomrule
      \end{tabular}
    \end{center}
  \label{tab:stats}
\end{table}

The above experiments involve high-dimensional parameter spaces with
as many as 1000 nuisance parameters. As a result, the assessment of
the statistics requiring bootstrapping is demanding in terms of
computational time and cost, even when parallel computing with a large
number of cores is used. For this reason, the number of bootstrap
samples is limited to $1000$ in all simulation experiments.

The three blocks of rows in Table~\ref{tab:tails} give the estimated tail
probabilities of the statistics of interest for the case $q=1000$ and
$m=8$ for all three models considered. This combination of $q$ and $m$ was
selected because it is the least extreme setting (compared to the most
extreme $q = 1000$, $m = 4$) where departures from the expected
behaviour in terms of the distribution of the statistics starts
becoming apparent; the results for all the other combinations of $q$
and $m$ are provided in the Supplementary Materials. The following
sections give a more detailed discussion on the figures shown in
Table~\ref{tab:tails}.

 \renewcommand{\baselinestretch}{1}
\begin{table}[th!]
  \caption{Empirical tail probabilities $\times 100$ for the
    statistics in Table~\ref{tab:stats} and all models considered in
    the simulation studies of Section~\ref{sec:simulations}. The
    figures shown have been rounded to 1 decimal and are for
    $q = 1000$ and $m = 8$.}
  \label{tab:tails}
  \begin{center}
    \begin{tabular}{llrrrrrr}
      \toprule
      & & \multicolumn{6}{c}{Nominal} \\
      \multicolumn{1}{c}{Model} & \multicolumn{1}{c}{Statistic} & $1.0$ & $2.5$ & $5.0$ & $95.0$ & $97.5$ & $99.0$ \\ \midrule
\multirow{6}{*}{Beta} &             $R$ &   0.0 &   0.0 &   0.0 &   0.0 &   0.0 &   0.0 \\
 &           $R^*$ &   0.7 &   1.8 &   3.8 &  93.7 &  96.8 &  98.8 \\
 &         $R^u$ &   0.8 &   1.9 &   4.1 &  94.0 &  97.0 &  98.7 \\
 &         $R^c$ &   1.0 &   2.3 &   4.8 &  95.0 &  97.4 &  99.1 \\
 &       $R^c_l$ &   1.1 &   2.5 &   5.1 &  94.7 &  97.3 &  98.9 \\
 & $R^c_{ls}$ &   0.9 &   2.3 &   4.8 &  95.1 &  97.5 &  99.0 \\ \midrule
 \multirow{6}{*}{Curved exponential family} &             $R$ & 100.0 & 100.0 & 100.0 & 100.0 & 100.0 & 100.0 \\
 &           $R^*$ &   1.4 &   3.5 &   6.9 &  96.6 &  98.3 &  99.4 \\
 &         $R^u$ &   0.6 &   1.8 &   4.0 &  95.0 &  97.7 &  99.2 \\
 &         $R^c$ &   1.2 &   3.3 &   6.4 &  96.2 &  98.2 &  99.4 \\
 &       $R^c_l$ &   1.5 &   3.6 &   7.1 &  95.8 &  98.0 &  99.2 \\
 & $R^c_{ls}$ &   1.3 &   3.2 &   6.5 &  96.3 &  98.2 &  99.4 \\ \midrule
\multirow{6}{*}{Truncated  regression} &             $R$ &   0.2 &   0.5 &   1.1 &  84.2 &  90.4 &  95.1 \\
 &           $R^*$ &   1.0 &   2.5 &   5.2 &  94.8 &  97.3 &  98.9 \\
 &         $R^u$ &   0.9 &   2.3 &   4.8 &  94.9 &  97.2 &  98.9 \\
 &         $R^c$ &   0.9 &   2.4 &   4.9 &  94.5 &  97.2 &  98.7 \\
 &       $R^c_l$ &   1.0 &   2.4 &   5.0 &  94.4 &  97.0 &  98.8 \\
 & $R^c_{ls}$ &   0.9 &   2.4 &   5.0 &  94.4 &  97.0 &  98.8 \\       
 \bottomrule
\end{tabular}
\end{center}
\end{table}
\renewcommand{\baselinestretch}{2}

\subsection{Beta model}
As a first example, we suppose that  $Y_{ij}$ has a beta distribution, with density function
\[
  g(y_{ij}; \mu_{i}, \phi) = \dfrac{1}{B\{\mu_{i}\phi, (1 -
    \mu_{i})\phi\}} y_{ij}^{\mu_{i}\phi - 1} (1 - y_{ij})^{(1 -
    \mu_{i})\phi - 1} \quad (0 < y_{ij} < 1) \,,
\]
where $B(\cdot)$ is the beta function. The parameter of interest is
$\psi = \log \phi$, whereas the stratum-specific nuisance parameters
are given by $\lambda_{i} = \log\{\mu_i / (1 - \mu_i)\}$. The
simulation experiments are carried out for $\psi_0 = \log(2)$ and the
elements of $\lambda_0$ are generated from a standard normal distribution
and  held fixed over all the replications.

\begin{figure}[t]
 \centering
 \includegraphics[width=1.1\textwidth]{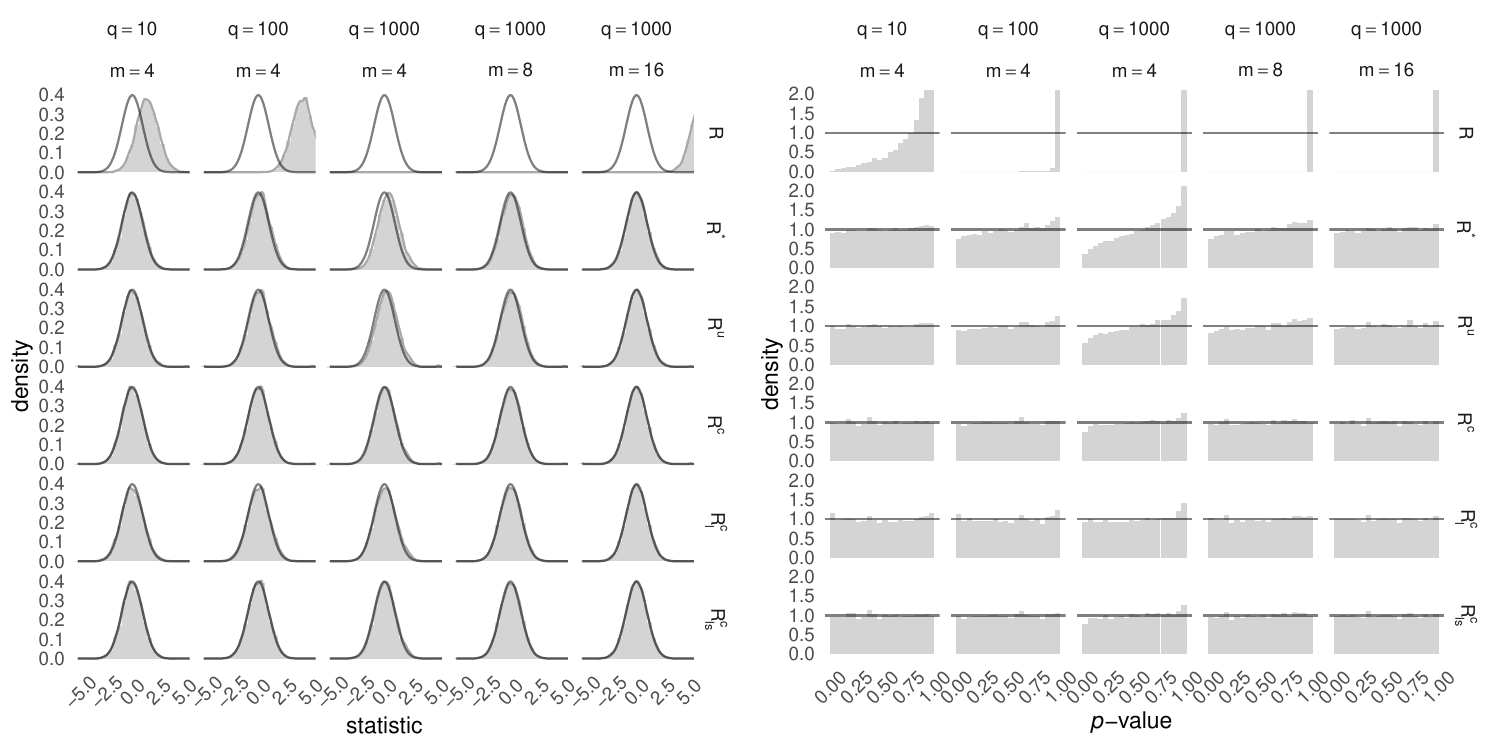}
 \caption{Beta model. Estimated null distribution of statistics (left)
   and estimated distribution of $p$-values (right) for the statistics
   in Table~\ref{tab:stats} for various combinations of $q$ and $m$.
   The ${\rm N}(0,1)$ and ${\rm Uniform}(0, 1)$ density functions are
   superimposed for statistics (left) and $p$-values (right).}
\label{fig:beta1_stat}
\end{figure}

The left panel of Figure~\ref{fig:beta1_stat} shows the empirical
densities for the statistics in
Table~\ref{tab:stats}. The performance of the statistics is evaluated
in terms of the closeness of their empirical density to the standard
normal density. This assessment is valid also for the constrained and
unconstrained bootstrap $p$-values, since they have been mapped into
the standard normal scale by the $\Phi^{-1}(\cdot)$ transformation.

The large location bias of the distribution of $R(\psi)$ is apparent
for all shown combinations of $q$ and $m$, and it becomes huge for
$q = 1000$ and $m \in \{4, 8\}$. All higher-order accurate statistics
result in a marked finite-sample correction, with $R^{*}(\psi)$ and
the unconstrained bootstrap illustrating some discrepancy from the
standard normal distribution
for large $q/m$ ratios, such as $q = 1000$  and $m \in \{4, 8\}$.
This is also apparent from the entries in Table~\ref{tab:tails}. 

From the right panel of Figure~\ref{fig:beta1_stat}, it is noticeable
that the $p$-values based on $R^{c}$, the location adjusted version
$R_{l}^{c}$ and the location-and-scale adjusted version $R_{ls}^{c}$
are all close to one another. Hence, the necessary adjustment for
making the distribution of $R(\psi)$ to be close to standard normal
appears to be mainly a location adjustment.

\subsection{Curved exponential family}
This example involves normally distributed random variables $Y_{ij}$,
each with mean $\exp(\lambda_{i})$ and variance
$\exp(\psi + \lambda_{i}/2)$. This model was studied in
\cite{sartori1999}, where it is pointed out that a marginal likelihood
for $\psi$ is not available. The simulation experiments are carried
out for $\psi_0 = \log(1/2)$ and the elements of $\lambda_0$ are
generated from a standard normal distribution and  held fixed over all the replications. 
 
\begin{figure}[t]
 \centering
 \includegraphics[width=1.1\textwidth]{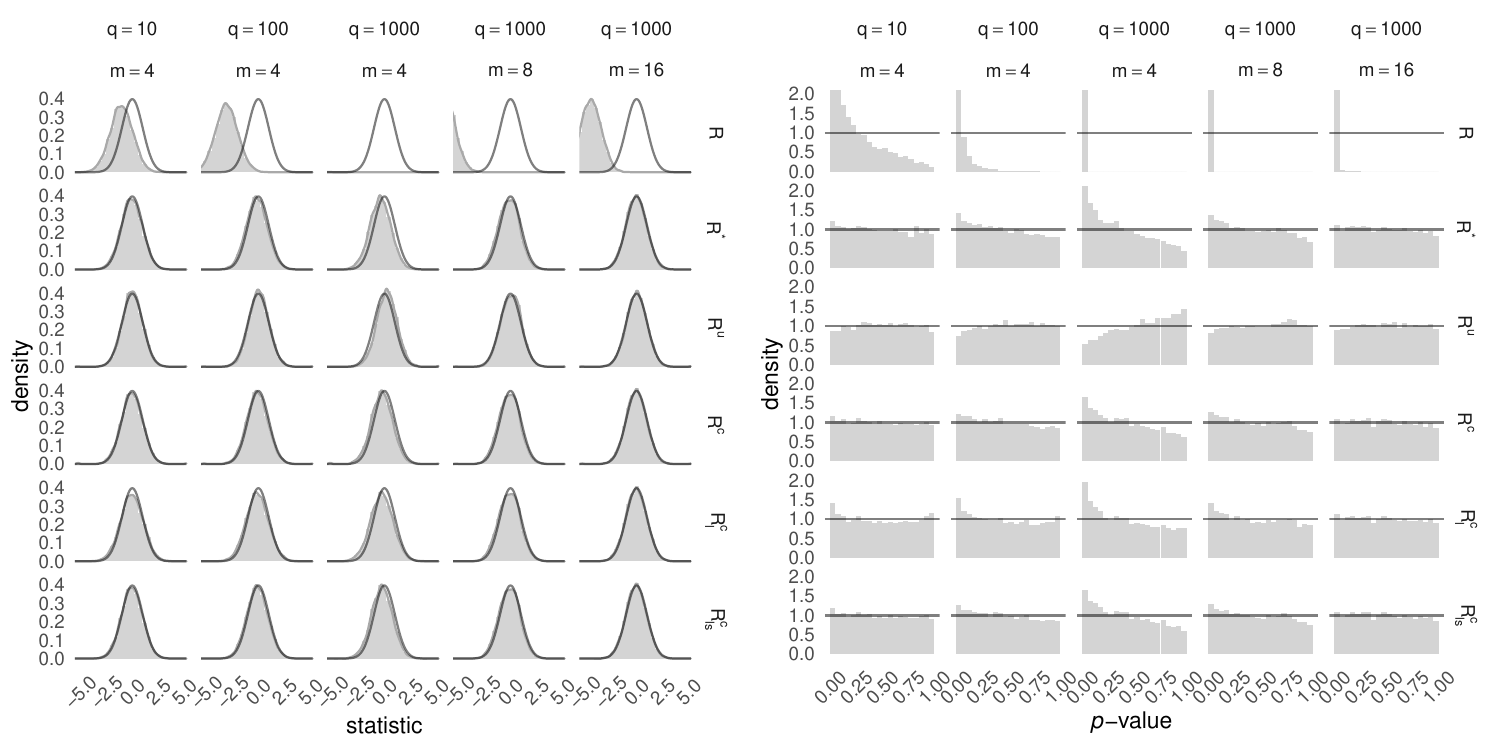}
 \caption{Curved exponential family model. Estimated null distribution
   of statistics (left) and estimated distribution of $p$-values
   (right) for the statistics in Table~\ref{tab:stats} for various
   combinations of $q$ and $m$. The ${\rm N}(0,1)$ and
   ${\rm Uniform}(0, 1)$ density functions are superimposed for
   statistics (left) and $p$-values (right).}
\label{fig:cef1_stat}
\end{figure}

The left panel in Figure \ref{fig:cef1_stat} shows the empirical
density functions of the statistics in Table~\ref{tab:stats}, and the
right panel shows the corresponding $p$-value distributions. 
Like in the previous example, 
the empirical, finite-sample distributions of $R(\psi)$ are far from
standard normal, while all the higher-order statistics perform
considerably better. The conclusions are similar to those from the
simulation experiments for the beta model, in that the
required adjustment to $R(\psi)$ seems to be a location
correction. The main difference is the fact that no statistic appears
to perform well for $(q, m) = (1000, 4)$; see also the empirical tail
probabilities in Table~\ref{tab:tails}.

\subsection{Truncated linear regression model}
The last example is taken from the econometric literature; see 
\cite{Greene:2004}, \cite{Bartolucci:2016} and the references therein.  We define the response variable 
$Y_{ij}$  to be distributed as $Y_{ij}^{*}$ conditionally on $y_{ij}^{*}>0$, with 
$$
y_{ij}^{*}=\lambda_{i } + x_{ij} \, \psi +  \varepsilon_{ij}\, , \quad i=1,\ldots,q, \,\,\,\,\,  j=1,\ldots,m\, ,
$$
where the error term $ \varepsilon_{ij}$ is standard normally distributed. For the simulation study, $\psi$ is set to  1
and the elements of $\lambda_0$ are generated from a standard normal distribution and  held fixed over all the replications. Likewise,  the 
values $x_{ij}$ are generated from a  standard normal distribution and  held fixed over all the replications.

\begin{figure}[t]
 \centering
 \includegraphics[width=1.1\textwidth]{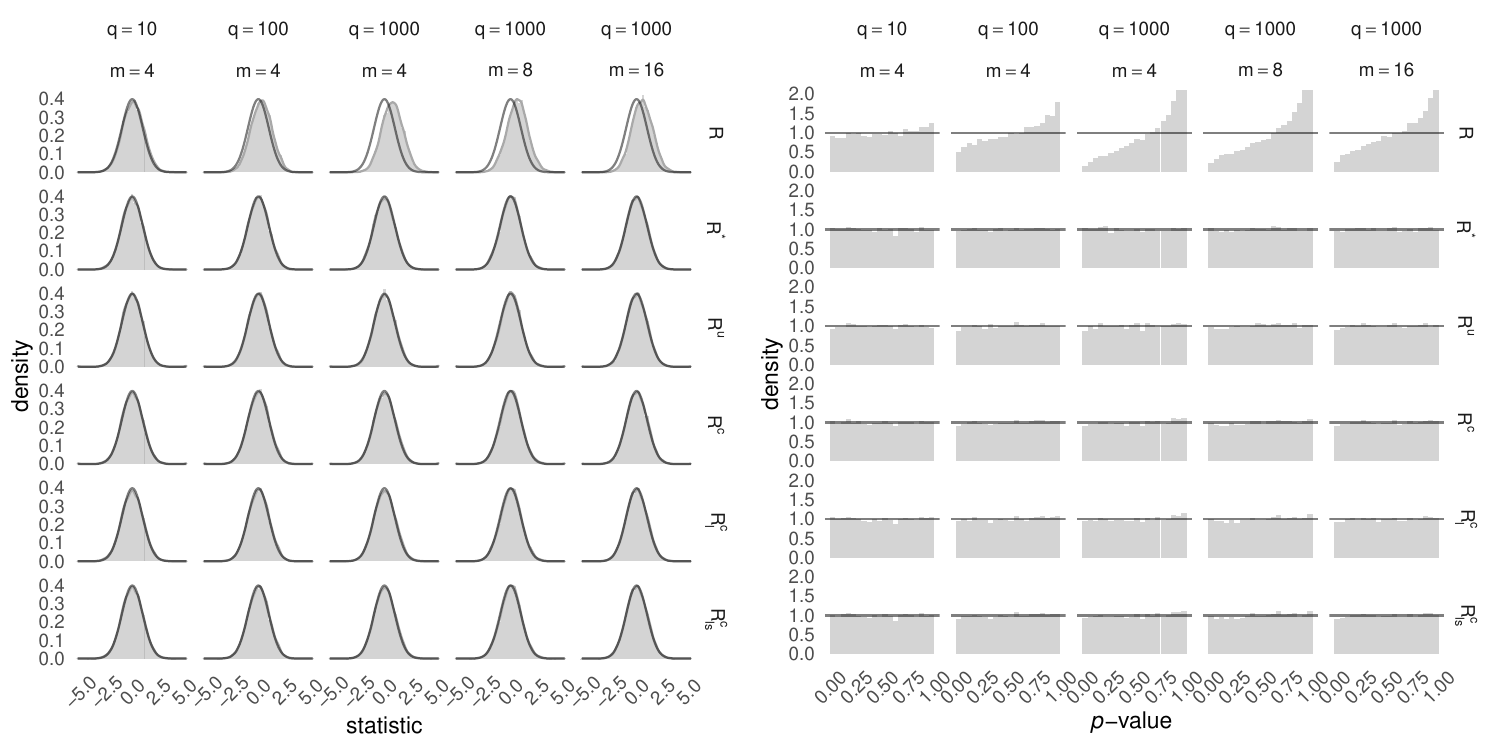}
 \caption{Truncated linear regression model. Estimated null distribution
   of statistics (left) and estimated distribution of $p$-values
   (right) for the statistics in Table~\ref{tab:stats} for various
   combinations of $q$ and $m$. The ${\rm N}(0,1)$ and
   ${\rm Uniform}(0, 1)$ density functions are superimposed for
   statistics (left) and $p$-values (right).}
\label{fig:treg_stat}
\end{figure}

The left panel in Figure \ref{fig:treg_stat} shows the empirical
density functions of the statistics in Table~\ref{tab:stats}, and the
right panel shows the corresponding $p$-value distributions.  Differently from the other examples, here the distribution of the
first-order statistics requires only a moderate adjustment even in the most extreme settings, and both the bootstrap-based
statistics as well the $R^{*}(\psi)$ statistic perform rather well, providing results very close to the target distributions.

\section{Concluding remarks}
\label{sec:discussion}
The contribution of this paper is to formally show that, in stratified
settings, inference based on either unconstrained or constrained
parametric bootstrap of usual likelihood pivots is effective in
recovering their inferential performance, even in rather extreme
settings, where the bias of the profile score renders vanilla
first-order inference invalid.

Unconstrained and constrained bootstrap for the signed likelihood
ratio root, the score statistic and the Wald statistic can both
recover inferential performance in stratified settings when
$q = O(m^\alpha)$, for $0 < \alpha < 3$. As in the case for
$\alpha = 0$ \citep{Lee.Young:2005}, when $0 < \alpha < 1$, constrained
bootstrap is seen to have a higher degree of asymptotic accuracy than
unconstrained bootstrap. On the other hand, the two bootstraps are
asymptotically equivalent when $1 \leq \alpha < 3$. The condition
$q = O(m^\alpha)$, for $0 < \alpha < 3$, is the same as the one found
in \cite{Sartori:2003} for validity of inference based on $R^*$ and on
the signed likelihood root computed from the modified profile
likelihood. 

The results of Section~\ref{sec:simulations} from the extensive
simulation studies for the finite-sample assessment of the performance
of constrained and unconstrained bootstrap are in par to what is
expected from theory. In extreme settings, like the beta model with
$(q, m) = (1000, 4)$, constrained bootstrap appears to perform
slightly better than unconstrained bootstrap. Furthermore, in all
simulation experiments we carried out and as $q/m$ diverges, the
inferential performance from constrained and unconstrained bootstrap
of first-order statistics seems to be deteriorating much slower than
that of $R^*$ and the signed likelihood root computed from the
modified profile likelihood (see also the Supplementary Materials). As
a result, the evidence from the simulation studies points out that
inference from parametric bootstrap is more resilient to increasing
$q/m$ than inference from well-used, analytically available
higher-order statistics, with constrained bootstrap being the most
accurate in extreme scenarios.

The theoretical developments in this paper do not immediately cover
situations where the random variables have discrete support, because
the Edgeworth expansion in~(\ref{Etrue}) can only be valid for models
with continuous support.  The impact of discreteness on the
performance of parametric bootstrap is examined in the Supplementary
Materials through a binomial matched pairs model. In particular, the
experimental setup of Section~\ref{sec:simulations} is used for a
stratified logistic regression model, where $Y_{ij}$ has a Bernoulli
distribution with probability
$\exp(\lambda_i + \psi x_{j})/\{1+\exp(\lambda_i + \psi x_{j})\}$,
with $x_{j} = 1$ for $j \in \{1, \ldots, m/2\}$ and $x_{j} = 0$ for
$j \in \{m/2 + 1, \ldots, m\}$. The results in Figures~S21-S24 and
Tables~S3-S11 in the Supplementary Materials indicate that the
equivalence between unconstrained and constrained bootstrap of the
first-order statistics in continuous models may not hold for discrete
settings. In those cases, despite that unconstrained bootstrap appears
to deliver a marked inferential improvement to first-order statistics,
constrained bootstrap, similarly to $R^*$, is found to perform
considerably better for most combinations of $q$ and $m$.

The simulation experiments in this paper have been carried out with
$1000$ bootstrap replications. This value is smaller than some of the
recommendations of millions of replications that have appeared in the
literature for standard asymptotics settings \citep{Young:2009,
  Diciccio.etal:2017}. For stratified settings with $\alpha > 1$ the
bootstrap adjustments have the role of recovering asymptotic
uniformity of $p$-values, rather that providing a small-sample
refinement of $p$-values that are asymptotically uniform. As a result,
use of a huge number of bootstrap replications is less essential, and
the few experiments we carried out with more than $1000$ bootstrap
replications are in support of that statement. More comprehensive
simulation studies to support that statement are unfortunately not
feasible with current computing capabilities.

%%%%%%%%%%%%%%%%%%%%%%%%%%%%%%%%%%%%%%%%%%%%%%%%%%%%%%%%%%%%%%%%%%%%%%%%%%%%%%%%%%%%%%%%%%%%%%%%%%%%%%%%%%%%%%%%%%%%%%%%%%%%
\section*{Supplementary Materials}
The Supplementary Materials provide the outputs from the simulation
experiments described in Section~\ref{sec:simulations}, for all models
and all combinations of statistics, $q$ and $m$. Outputs are also
provided for other models, given by a gamma model, a Behrens-Fisher model and 
by
the logistic regression model described in
Section~\ref{sec:discussion}. The outputs include null distributions
of the various statistics and distributions of $p$-values, through
extended versions of the
Figures~\ref{fig:beta1_stat}-\ref{fig:treg_stat}, and empirical tail
probabilities, through extended versions of Table~\ref{tab:tails}.

%%%%%%%%%%%%%%%%%%%%%%%%%%%%%%%%%%%%%%%%%%%%%%%%%%%%%%%%%%%%%%%%%%%%%%%%%%%%%%%%%%%%%%%%%%%%%%%%%%%%%%%%%%%%%%%%%%%%%%%%%%%%
\section*{Acknowledgements}
We thank two anonymous Referees and the Associate Editor for their constructive comments which enabled us  to significantly improve our work. 
Ioannis Kosmidis has been supported by The Alan Turing Institute under
the EPSRC grant EP/N510129/1. Part of this work was completed during a
visit of Ruggero Bellio, Alessandra Salvan and Nicola Sartori to The
Alan Turing Institute. Alessandra Salvan and Nicola Sartori have been
supported by the University of Padova under grants BIRD185955 and BIRD203991.

%%%%%%%%%%%%%%%%%%%%%%%%%%%%%%%%%%%%%%%%%%%%%%%%%%%%%%%%%%%%%%%%%%%%%%%%%%%%%%%%%%%%%%%%%%

\setcounter{equation}{0}
\renewcommand{\theequation}{A\arabic{equation}}

\section*{Appendix}
\subsection*{Asymptotic orders in (\ref{ns2003})}

\noindent
The following representation from \citet[][Appendix]{Sartori:2003}
will be used to determine the order of quantities in a stratified
setting. Let $\mu_i$ and $\sigma^2_i$  denote mean and variance of
independent the random variables $X_1, \ldots, X_q$. Then
\begin{equation}\label{CLT}
\sumq X_i=O_p\left(\sumq \mu_i\right)+O_p\left(\sqrt{\sumq\sigma^2_i}\right)\,.
\end{equation}
We have $U_\psi=\sum_{i=1}^q U^i_\psi$, where $U_\psi^i$ is the
contribution to $U_\psi$ from the $i$th stratum, and
$U_\lambda=(U_{\lambda_1},\ldots, U_{\lambda_q})^\top$.  Here and in
the following, when the argument is omitted, evaluation at $\theta$ is
understood.

The terms on the right-hand side of (\ref{ns2003}) are seen to be of order $O_p(m^{(\alpha+1)/2})$, 
$O_p(m^{\alpha})$ and $O_p(m^{\alpha-1})$, respectively. Indeed,  
 using (\ref{CLT}), we have $U_{\psi|\lambda}=  \sum_{i=1}^q U_{\psi|\lambda_i}=O_p(m^{(\alpha+1)/2})$, with $U_{\psi|\lambda_i}= U_\psi^i - i_{\psi\lambda_i} i_{\lambda_i\lambda_i}^{-1} U_{\lambda_i}$ being $\E_\theta(U_{\psi|\lambda_i})=0$ and $\V_\theta(U_{\psi|\lambda_i})=i_{\psi\psi |\lambda_i}=O(m)$. Note that $i_{\psi\psi |\lambda}= \V_\theta(U_{\psi|\lambda})=\sumq i_{\psi\psi |\lambda_i}$. Similarly, we have
$B=\sumq B^i(\psi,\lambda_i)=O_p(m^\alpha)$, where $B^i(\psi,\lambda_i)$ is the term of order $O_p(1)$ of the expansion of the profile score in the $i$th stratum, having both mean and variance of order $O(1)$. The same additivity property holds for $b(\theta)$, so that $b(\theta)=\sumq b^i(\psi,\lambda_i)=O(m^\alpha)$. 
 Finally, the remainder term is $Re=\sumq Re^i(\psi,\lambda_i)$, with $Re^i(\psi,\lambda_i)$ having mean and variance of order $O(m^{-1})$, so that $Re=O_p(m^{\max\{\alpha-1, (\alpha-1)/2\}})=O_p(m^{\alpha-1})$ when $\alpha>1$.

\subsection*{Derivation of (\ref{ESexp}) and (\ref{VSexp})}

As a first step, consider the expansion
\begin{equation}\label{ippl}
i_{\psi\psi |\lambda}(\hat\theta_\psi)=i_{\psi\psi |\lambda}+C + O_p(m^{\alpha-1})\,,
\end{equation}
with
\begin{equation*}
C= \sumq \frac{d}{d\lambda_i}i_{\psi\psi |\lambda_i}(\hat\lambda_{i\psi}-\lambda_i)+
\frac{1}{2} \sumq  \frac{d^2}{d\lambda_i^2}i_{\psi\psi |\lambda_i}(\hat\lambda_{i\psi}-\lambda_i)^2= O_p(m^\alpha)\,,
\end{equation*}
where when $\alpha>1$ both terms in $C$ are of the same order, which is again determined using (\ref{CLT}).
Hence,
\begin{equation}\label{p_info_exp}
\{i_{\psi\psi |\lambda}(\hat\theta_\psi)\}^{-1/2} =
i_{\psi\psi |\lambda}^{-1/2}
 \left\{ 1 -\frac{1}{2} \frac{C}{i_{\psi\psi |\lambda}}+O_p(m^{-2}) \right\}\,,
\end{equation}
with $C/i_{\psi\psi |\lambda}=O_p(m^{-1})$.

Using (\ref{ns2003}) and (\ref{p_info_exp}), 
\begin{eqnarray}
S(\psi) &=& i_{\psi\psi |\lambda}^{-1/2} \left\{ U_{\psi|\lambda}+B  +Re\right\} \left\{1 -\frac{1}{2} \frac{C }{i_{\psi\psi |\lambda }}+O_p(m^{-2})\right\}\nonumber \\
&=& \frac{U_{\psi|\lambda}}{i_{\psi\psi |\lambda}^{1/2}} +
\frac{B }{i_{\psi\psi |\lambda}^{1/2}} +\frac{Re}{i_{\psi\psi |\lambda}^{1/2}} 
-\frac{1}{2} \frac{U_{\psi|\lambda}C }{i_{\psi\psi |\lambda}^{3/2}} 
-\frac{1}{2} \frac{B\, C }{i_{\psi\psi |\lambda}^{3/2}}  -\frac{1}{2} \frac{Re\,C }{i_{\psi\psi |\lambda}^{3/2}} \nonumber \\
&&+ O_p(m^{-2})+O_p(m^{(\alpha-5)/2})+O_p(m^{(\alpha-7)/2})\,,  \label{Sexp}
\end{eqnarray}
where $O_p(m^{-2})+O_p(m^{(\alpha-5)/2})+O_p(m^{(\alpha-7)/2})= O_p(m^{(\alpha-5)/2})$ as long as $\alpha>1$.
The term of order $O_p(m^{(\alpha-5)/2})$ is given by 
$
i_{\psi\psi |\lambda}^{-1/2}B$ times the term of order $O_p(m^{-2})$ in (\ref{p_info_exp}). 
Its expectation is of order $O(m^{(\alpha-5)/2})$.
The orders of   terms in (\ref{Sexp}) are as follows:
\begin{eqnarray*}
&& \frac{U_{\psi|\lambda}}{i_{\psi\psi |\lambda}^{1/2}} = O_p(1)\,, \quad
\frac{B }{i_{\psi\psi |\lambda}^{1/2}}= O_p(m^{(\alpha-1)/2})\,, \quad 
\frac{Re}{i_{\psi\psi |\lambda}^{1/2}} = O_p(m^{(\alpha-3)/2})\,,\\
&&\frac{1}{2} \frac{U_{\psi|\lambda}C }{i_{\psi\psi |\lambda}^{3/2}} = O_p(m^{-1})=o_p(1) \,,\quad
\frac{1}{2} \frac{B\,C}{i_{\psi\psi |\lambda}^{3/2}} = O_p(m^{(\alpha-3)/2})\,,\\
&&
-\frac{1}{2} \frac{Re\,C}{i_{\psi\psi |\lambda}^{3/2}} = O_p(m^{(\alpha-5)/2})\,.
\end{eqnarray*}

Expansion (\ref{ESexp}) for $\E_\theta(S(\psi))$ is obtained using
(\ref{Sexp}) and recalling that $b(\theta)=O(m^\alpha)$. We have
\begin{eqnarray*}
&& \E_\theta\left(\frac{U_{\psi|\lambda}}{i_{\psi\psi |\lambda}^{1/2}}\right) = 0\,, \quad
 \E_\theta\left(\frac{B }{i_{\psi\psi |\lambda}^{1/2}}\right)=\frac{b(\theta)}{i_{\psi\psi |\lambda}^{1/2}}= O(m^{(\alpha-1)/2})\,, \\ 
&& \E_\theta\left(\frac{Re}{i_{\psi\psi |\lambda}^{1/2}}\right) = O(m^{(\alpha-3)/2})\,,\quad
 \E_\theta\left(\frac{1}{2} \frac{U_{\psi|\lambda}C }{i_{\psi\psi |\lambda}^{3/2}}\right) = O(m^{-(\alpha+3)/2})=o(1)\,, \\
&& \E_\theta\left(\frac{1}{2} \frac{B \,C }{i_{\psi\psi |\lambda}^{3/2}}\right) = O_p(m^{(\alpha-3)/2})\,,\quad
\E_\theta\left(-\frac{1}{2} \frac{Re\,C }{i_{\psi\psi |\lambda}^{3/2}}\right) = O(m^{(\alpha-5)/2})\,,
\end{eqnarray*}
giving (\ref{ESexp}) with
\begin{equation}\label{M1}
M_1(\theta)= \E_\theta\left(\frac{Re}{i_{\psi\psi |\lambda}^{1/2}}\right) + 
\E_\theta\left(\frac{1}{2} \frac{B\,C}{i_{\psi\psi |\lambda}^{3/2}}\right) =O(m^{(\alpha-3)/2})\,.
\end{equation}

Expansion (\ref{VSexp}) for $\V_\theta(S(\psi))$ is also obtained
using (\ref{Sexp}).  In particular, the leading term has variance
equal to 1, and, using a standard expansion for the stratum profile
score $U_\psi^i(\psi,\hat\lambda_{i\psi})$ \citep[see e.g][formula
(8.88)]{Pace.Salvan:1997}, $\Cov_\theta (U_{\psi|\lambda}, B)$ and
$\V_\theta(B)$ are easily seen to be of order $O(m^\alpha)$.  Further terms of (\ref{Sexp}) give
contributions to the variance of order $O(m^{-2})$.
   
Higher order cumulants of $S(\psi)$, $r=3,4,\ldots$, have the form
\[
\kappa_r(S(\psi))=\frac{O(m^{\alpha+1})}{O(m^{r(\alpha+1)/2})}= O(m^{(\alpha+1)(1-r/2)})= O(n^{1-r/2})
\]
as in standard asymptotics.

\subsection*{Derivation of (\ref{BS}) and (\ref{v})}

\noindent
Let $\bRe=  \E_\theta(Re)$ and $\bBC=\E_\theta(B\,C)$. Then, from (\ref{ESexp}) and (\ref{M1}),  
\[
M(\hat\theta_\psi)=\left\{ i_{\psi\psi |\lambda}(\hat\theta_\psi)\right\}^{-1/2}
\left\{ 
b(\hat\theta_\psi) + \bRe(\hat\theta_\psi) +\frac{1}{2}\frac{1}{ i_{\psi\psi |\lambda}(\hat\theta_\psi) }\bBC(\hat\theta_\psi)
\right\} +O_p(m^{(\alpha-5)/2})\,,
\]
where $\bRe(\hat\theta_\psi)$ and
$ i_{\psi\psi |\lambda}(\hat\theta_\psi)^{-1}\bBC(\hat\theta_\psi)$ are of order $O(m^{\alpha-1})$.
Now,
\begin{equation}\label{expBP}
b(\hat\theta_\psi)=b(\theta)+b_1(\theta)+O_p(m^{\alpha-2})\,,
\end{equation}
where 
\begin{equation}\label{b1}
b_1(\theta)=\sumq b^i_{\lambda_i}(\psi,\lambda_i)(\hat\lambda_{i\psi}-\lambda_i)+
\frac{1}{2} \sumq b^i_{\lambda_i\lambda_i}(\psi,\lambda_i)(\hat\lambda_{i\psi}-\lambda_i)^2\,,
\end{equation}
and $b^i_{\lambda_i}(\psi,\lambda_i)=\partial b^i(\psi,\lambda_i)/\partial\lambda_i$, and so on. 
Using (\ref{CLT}), and being $b^i_{\lambda_i}(\psi,\lambda_i)$
 and $b^i_{\lambda_i\lambda_i}(\psi,\lambda_i)$ both of order $O(1)$,  
\begin{equation*}\label{BS1}
\sumq b^i_{\lambda_i}(\psi,\lambda_i)(\hat\lambda_{i\psi}-\lambda_i) =O_p(m^{\alpha-1})+
O(m^{(\alpha-1)/2}) 
\end{equation*}
and 
\begin{equation*}\label{BS2}
\sumq b^i_{\lambda_i\lambda_i}(\psi,\lambda_i)(\hat\lambda_{i\psi}-\lambda_i)^2
=O_p(m^{\alpha-1})+O_p(m^{(\alpha-2)/2})\,.
\end{equation*}
The remainder in (\ref{expBP}) is of order $O_p(m^{\alpha-2})+O_p(m^{(\alpha-3)/2})=O_p(m^{\alpha-2})$, when $\alpha>1$.
Moreover, $\bRe(\hat\theta_\psi)=\bRe  + O_p(m^{\alpha-2})$ and
$i_{\psi\psi |\lambda}(\hat\theta_\psi)^{-1}\bBC(\hat\theta_\psi) =  i_{\psi\psi |\lambda}^{-1}\bBC + O_p(m^{\alpha-2})$.

Using (\ref{p_info_exp}),  we get
\begin{equation}\label{Mconstr}
M(\hat\theta_\psi)=
i_{\psi\psi |\lambda}^{-1/2} b(\theta)
+  \tilde{M}_1
+O_p\left(  m^{-\min\{1, (5-\alpha)/2\}}\right)\,, 
\end{equation}
with
\[
\tilde{M}_1= i_{\psi\psi |\lambda}^{-1/2}\left\{
b_1(\theta)-\frac{C\,b(\theta)}{2i_{\psi\psi |\lambda}}+\bRe +\frac{\bBC}{2 i_{\psi\psi |\lambda}}
\right\}\,,
\]
which is of order $O_p(m^{(\alpha-3)/2})$ because all terms are of the
same order.

Therefore, (\ref{Mconstr}), (\ref{ESexp}) and (\ref{M1}) give  (\ref{BS}) with
\begin{equation}\label{Delta}
\Delta=\tilde{M}_1-M_1(\theta)
=\frac{b_1(\theta)}{i_{\psi\psi |\lambda}^{1/2}}  
-  \frac{C\,b(\theta)}{2\,i_{\psi\psi |\lambda}^{3/2}}  
\end{equation}
that is of order $O_p(m^{(\alpha-3)/2})$.

To obtain expansion (\ref{v}) recall that in (\ref{VSexp})
\[
v(\theta)=(\V_\theta(B) + 2\,\E_\theta(U_{\psi|\lambda} B))/i_{\psi\psi |\lambda}\,.
\]
The numerator of $v(\hat\theta_\psi)$ is equal to $\V_\theta(B) + 2\,\E_\theta(U_{\psi|\lambda} B)$ plus a term of order $O_p(m^{\alpha-1})$. From (\ref{ippl}),
\begin{equation*}\label{ippl1}
\frac{1}{i_{\psi\psi |\lambda}(\hat\theta_\psi)}=\frac{1}{i_{\psi\psi |\lambda}}\left\{ 1+O_p(m^{-1})  \right\}\,.
\end{equation*}
which gives  (\ref{v}).

\subsection*{Derivation of (\ref{res2})}

First, from (\ref{E1}), we have 
\begin{equation*}\label{E2}
F_{\hat\theta}(x)= \Phi\left(x^*(\hat\theta)\right)+ O_p\left( m^{-\min\left(2, \frac{\alpha+1}{2}\right)} \right)\,.
\end{equation*}
In order to obtain expansions of $M(\hat\theta)$ and $v(\hat\theta)$
around $\theta$, we use the fact that, when $\alpha>1$,
$\hat\psi-\psi=O_p(m^{-1})$ \citep{Sartori:2003}.  This implies that
an expansion for $F_{\hat\theta}(x)$ of the form (\ref{CBE}) holds
with a different $\Delta$ term, which is still of order
$O_p(m^{(\alpha-3)/2})$.

In order to obtain an expansion for $M(\hat\theta)$ we follow the same
steps as in (\ref{expBP})--(\ref{Delta}), giving (\ref{BS}).  In
particular, we have
\begin{equation*}\label{expBPconv}
b(\hat\theta)= b(\theta)+b_{2}(\theta)+O_p(m^{\alpha-2}) \, ,
\end{equation*}
where
\begin{eqnarray}
b_2(\theta)=b_2(\psi,\lambda)&=&\sumq b^i_{\psi} \, (\hat\psi-\psi)+\sumq b^i_{\lambda_i}(\hat\lambda_i-\lambda_i) 
+\frac{1}{2} \sumq b^i_{\lambda_i\lambda_i} (\hat\lambda_i-\lambda_i)^2
\nonumber \\
&& +\frac{1}{2} \sumq b^i_{\psi\psi}(\hat\psi-\psi)^2+\sumq b^i_{\psi\lambda_i} (\hat\lambda_i-\lambda_i) (\hat\psi-\psi)\,. \label{expBPconv}
\end{eqnarray}
From \citet[][below formula (9)]{Sartori:2003}, with $\alpha>1$,
$\hat\psi-\psi=O_p(m^{-1})$, so that the first three summands on the
right hand side of the last formula are of order $O_p(m^{\alpha-1})$,
while the remaining two are of order $O_p(m^{\alpha-2})$. This leads
to
 \begin{equation}\label{BShat}
M(\hat\theta)= M(\theta)+\Delta_1+O_p\left(m^{-\min\left(1, \frac{5-\alpha}{2}\right)}\right) \,,
\end{equation}
where the term $\Delta_1$ is of order $O_p(m^{(\alpha-3)/2})$, as its
expected value, because the leading terms in (\ref{expBPconv}) are of
the same order as $b_1(\theta)$ in (\ref{expBP}).

Using (\ref{BShat}) and an expansion similar to (\ref{v}) we obtain
\[
x^*(\hat\theta)=
x^*(\theta)+O_p(m^{(\alpha-3)/2})\,,
\]
so that the same error as in (\ref{CBE}) holds also for unconstrained bootstrap, i.e.
\begin{equation}\label{BE}
F_{\hat\theta}(x)=
F_{\theta}(x)-\phi(x^*(\theta)) \Delta_1+ O_p(m^{-1})\,.
\end{equation}
The steps leading from (\ref{BE}) to (\ref{res2}) are the same as
those from (\ref{CBE}) to (\ref{res1}).

\bibliographystyle{chicago}
\bibliography{prepivoting}

\begin{thebibliography}{}

\bibitem[\protect\citeauthoryear{Barndorff-Nielsen}{Barndorff-Nielsen}{1986}]{Barndorff-Nielsen:1986}
Barndorff-Nielsen, O.~E. (1986).
\newblock Inference on full or partial parameters based on the standardized
  signed log likelihood ratio.
\newblock {\em Biometrika\/}~{\em 73}, 307--322.

\bibitem[\protect\citeauthoryear{Bartolucci, Bellio, Salvan, and
  Sartori}{Bartolucci et~al.}{2016}]{Bartolucci:2016}
Bartolucci, F., R.~Bellio, A.~Salvan, and N.~Sartori (2016).
\newblock Modified profile likelihood for fixed-effects panel data models.
\newblock {\em Econometric Reviews\/}~{\em 35}, 1271--1289.

\bibitem[\protect\citeauthoryear{DiCiccio, Kuffner, and Young}{DiCiccio
  et~al.}{2017}]{Diciccio.etal:2017}
DiCiccio, T.~J., T.~A. Kuffner, and G.~A. Young (2017).
\newblock The formal relationship between analytic and bootstrap approaches to
  parametric inference.
\newblock {\em Journal of Statistical Planning and Inference\/}~{\em 191},
  81--87.

\bibitem[\protect\citeauthoryear{DiCiccio, Martin, and Stern}{DiCiccio
  et~al.}{2001}]{DiCiccio.etal:2001}
DiCiccio, T.~J., M.~A. Martin, and S.~E. Stern (2001).
\newblock Simple and accurate one-sided inference from signed roots of
  likelihood ratios.
\newblock {\em The Canadian Journal of Statistics\/}~{\em 29}, 67--76.

\bibitem[\protect\citeauthoryear{Greene}{Greene}{2004}]{Greene:2004}
Greene, W. (2004).
\newblock The behaviour of the maximum likelihood estimator of limited
  dependent variable models in the presence of fixed effects.
\newblock {\em The Econometrics Journal\/}~{\em 7}, 98--119.

\bibitem[\protect\citeauthoryear{Lee and Young}{Lee and
  Young}{2005}]{Lee.Young:2005}
Lee, S. M.~S. and G.~A. Young (2005).
\newblock Parametric bootstrapping with nuisance parameters.
\newblock {\em Statistics and Probability Letters\/}~{\em 71}, 143--153.

\bibitem[\protect\citeauthoryear{Pace and Salvan}{Pace and
  Salvan}{1997}]{Pace.Salvan:1997}
Pace, L. and A.~Salvan (1997).
\newblock {\em Principles of Statistical Inference from a Neo-Fisherian
  Perspective}.
\newblock Singapore: World Scientific.

\bibitem[\protect\citeauthoryear{Pierce and Bellio}{Pierce and
  Bellio}{2017}]{Pierce.Bellio:2017}
Pierce, D.~A. and R.~Bellio (2017).
\newblock Modern likelihood-frequentist inference.
\newblock {\em International Statistical Review\/}~{\em 85}, 519--541.

\bibitem[\protect\citeauthoryear{Sartori}{Sartori}{2003}]{Sartori:2003}
Sartori, N. (2003).
\newblock Modified profile likelihoods in models with stratum nuisance
  parameters.
\newblock {\em Biometrika\/}~{\em 90}, 533--549.

\bibitem[\protect\citeauthoryear{Sartori, Bellio, Salvan, and Pace}{Sartori
  et~al.}{1999}]{sartori1999}
Sartori, N., R.~Bellio, A.~Salvan, and L.~Pace (1999).
\newblock The directed modified profile likelihood in models with many nuisance
  parameters.
\newblock {\em Biometrika\/}~{\em 86}, 735--742.

\bibitem[\protect\citeauthoryear{Severini}{Severini}{2000}]{Severini:2000}
Severini, T.~A. (2000).
\newblock {\em Likelihood Methods in Statistics}.
\newblock Oxford: Oxford University Press.

\bibitem[\protect\citeauthoryear{Stern}{Stern}{2006}]{Stern:2006}
Stern, S.~E. (2006).
\newblock {Simple and accurate one-sided inference based on a class of
  M-estimators}.
\newblock {\em Biometrika\/}~{\em 93}, 973--987.

\bibitem[\protect\citeauthoryear{Young}{Young}{2009}]{Young:2009}
Young, G.~A. (2009).
\newblock Routes to higher-order accuracy in parametric inference.
\newblock {\em Australian {\&} New Zealand Journal of Statistics\/}~{\em 51},
  115--126.

\bibitem[\protect\citeauthoryear{Young and Smith}{Young and
  Smith}{2005}]{Young.Smith:2005}
Young, G.~A. and R.~L. Smith (2005).
\newblock {\em Essentials of Statistical Inference}.
\newblock Cambridge: Cambridge University Press.

\end{thebibliography}



\section*{Supplementary material}

where

$$b_2(\theta) = b_2(\psi, \lambda) = \sum_{i=1}^{q} b_{\psi}^{i} (\hat{\psi} - \psi) + \sum_{i=1}^{q} b_{\lambda_i}^{i}(\hat{\lambda}_i - \lambda_i) + \frac{1}{2}\sum_{i=1}^{q} b_{\lambda_i \lambda_i}^{i}(\hat{\lambda}_i - \lambda_i)^2$$
$$+ \frac{1}{2}\sum_{i=1}^{q} b_{\psi\psi}^{i}(\hat{\psi} - \psi)^2 + \sum_{i=1}^{q} b_{\psi\lambda_i}^{i}(\hat{\lambda}_i - \lambda_i)(\hat{\psi} - \psi). \quad (S1.23)$$

From Sartori (2003, below formula (9)), with $\alpha < 1$, $\hat{\psi} - \psi = O_p(m^{-(\alpha+1)/2})$, so that leading terms of $b_2(\theta)$ are the first two summands on the right-hand side of (S1.23), which are of order $O_p(m^{(\alpha-1)/2})$. The other terms give a contribution of order $O_p(m^{\alpha-1})$. This leads to

$$M(\hat{\theta}) = M(\theta) + \Delta_1 + O_p(m^{(\alpha-3)/2}), \quad (S1.24)$$

where the term $\Delta_1$ has form similar to (S1.14) and is of order $O_p(m^{-1})$, with expectation of order $O(m^{(\alpha-3)/2})$, because the leading term in (S1.23) are of the same order as $b_1(\theta)$ in (S1.16).

The remaining steps are analogous to those from (S1.20) with $\Delta$ replaced by $\Delta_1$.

## S2 Full set of simulation results

We report the full set of results of the simulation studies described in the main text. In addition to the models described there, which include also the logistic regression model described in the final section, two further models are considered, given by a gamma model and a Behrens-Fisher model, respectively. Some details about them are given as follows.

**Gamma model**

We take a stratified gamma model with common shape parameter, for which $Y_{ij}$ has density function

$$g(y_{ij}; \alpha, \beta_i) = \frac{1}{\Gamma(\alpha)\, \beta_i^{\alpha}}\, y_{ij}^{\alpha-1}\, \exp\left\{-\frac{y_{ij}}{\beta_i}\right\} \quad (y_{ij} > 0),$$

where $\Gamma(\cdot)$ is the gamma function. The parameter of interest is $\psi = \log(\alpha)$ and the stratum-specific nuisance parameters are given by $\lambda_i = \log(\beta_i)$. The simulation data sets where generated for $\psi_0 = \log(2)$ and the elements of $\lambda_0$ are fixed to the logarithm of random draws from an exponential distribution with rate $1/2$ and held fixed over all the replications.

The results for this model, reported in the following, are similar to those for the beta model.

**Behrens-Fisher model**

Suppose that $Y_{ij}$ are normally distributed, with common mean $\psi$ across all strata, and stratum-specific variances $\exp(\lambda_i)$. This example was studied in Young (2009), where it is illustrated that $R^*(\psi)$ is performing worse than constrained bootstrapping of $R(\psi)$ in a simulation study with $q = 20$ strata and the stratum sample size varying between 3 and 10. The simulation experiments are carried out for $\psi_0 = 0$ and $\lambda_0$ generated from a uniform distribution in $(0, 1)$ and held fixed over all the replications.

The findings in Young (2009) are confirmed also in the more extensive simulation studies here, whose results are reported in the following. For the Behrens-Fisher model, the location-adjustment to $R(\psi)$ does not have much of an effect to the distribution of the statistic. On the other hand, the location-and-scale adjusted statistic recovers inferential performance. So, it seems that only a scale adjustment is sufficient to recover first-order inferential performance. Indeed, for this model the profile score bias is exactly equal to zero, so that the distribution of the first-order statistic $R(\psi)$ is essentially centred around zero. Also, both constrained and unconstrained bootstrap perform extremely well across all the simulation scenarios, and, along with $R_{ls}^c$, they outperform $R^*(\psi)$ in the most extreme settings.

A remarkable observation in this model is the excellent accuracy of the score statistic $S(\psi)$ for all combination of $q$ and $m$ considered. This can be explained by the fact that the leading term $U_{\psi|\lambda}$ of $U_p$ in (S1.3) is exactly normally distributed. Moreover, $M(\theta)$ and $v(\theta)$, in (S1.4) and (S1.6) respectively, are equal to zero.

## Results

For each model, simulation results for the 23 statistics are visualized. For plotting purposes, the statistics have been split in two groups. Table 1 lists the first group, and Table 2 the second one.

Table 1: Statistics computed in the simulation experiments (first group for plotting).

| Statistic | Plotting Symbol | Description |
|---|---|---|
| $R(\psi)$ | $R$ | Signed likelihood root |
| $S(\psi)$ | $S$ | Score statistic |
| $T(\psi)$ | $T$ | Wald statistic |
| $\Phi^{-1}\{\hat{p}_1^T(\psi)\}$ | $T^u$ | Transformed $p$-value from unconstrained bootstrap of $T(\psi)$ |
| $\Phi^{-1}\{\hat{p}_2^T(\psi)\}$ | $T^c$ | Transformed $p$-value from constrained bootstrap of $T(\psi)$ |
| $\Phi^{-1}\{\hat{p}_1^S(\psi)\}$ | $S^u$ | Transformed $p$-value from unconstrained bootstrap of $S(\psi)$ |
| $\Phi^{-1}\{\hat{p}_2^S(\psi)\}$ | $S^c$ | Transformed $p$-value from constrained bootstrap of $S(\psi)$ |
| $\Phi^{-1}\{\hat{p}_1^R(\psi)\}$ | $R^u$ | Transformed $p$-value from unconstrained bootstrap of $R(\psi)$ |
| $\Phi^{-1}\{\hat{p}_2^R(\psi)\}$ | $R^c$ | Transformed $p$-value from constrained bootstrap of $R(\psi)$ |
| $R^*(\psi)$ | $R^*$ | Modified signed likelihood root |
| $R_m(\psi)$ | $R_m$ | Signed likelihood root computed from the modified profile likelihood |

Table 2: Statistics computed in the simulation experiments (second group for plotting). The mean $\hat{\mu}^R$ and the standard deviation $\hat{\sigma}^R$ of $R(\psi)$ are estimated through unconstrained bootstrap, by simulating from the model at $\theta = \hat{\theta}$, and likewise for $S(\psi)$ and $T(\psi)$. The mean $\tilde{\mu}^R$ and the standard deviation $\tilde{\sigma}^R$ of $R(\psi)$ are estimated through constrained bootstrap, by simulating from the model at $\theta = \hat{\theta}_\psi$, and likewise for $S(\psi)$ and $T(\psi)$.

| Statistic | Plotting Symbol | Description |
|---|---|---|
| $T(\psi) - \hat{\mu}^T$ | $T_l^u$ | Location adjusted $T(\psi)$, unconstrained bootstrap |
| $T(\psi) - \tilde{\mu}^T$ | $T_l^c$ | Location adjusted $T(\psi)$, constrained bootstrap |
| $S(\psi) - \hat{\mu}^S$ | $S_l^u$ | Location adjusted $S(\psi)$, unconstrained bootstrap |
| $S(\psi) - \tilde{\mu}^S$ | $S_l^c$ | Location adjusted $S(\psi)$, constrained bootstrap |
| $R(\psi) - \hat{\mu}^R$ | $R_l^u$ | Location adjusted $R(\psi)$, unconstrained bootstrap |
| $R(\psi) - \tilde{\mu}^R$ | $R_l^c$ | Location adjusted $R(\psi)$, constrained bootstrap |
| $(T(\psi) - \hat{\mu}^T)/\hat{\sigma}^T$ | $T_{ls}^u$ | Location-and-scale adjusted $T(\psi)$, unconstrained bootstrap |
| $(T(\psi) - \tilde{\mu}^T)/\tilde{\sigma}^T$ | $T_{ls}^c$ | Location-and-scale adjusted $T(\psi)$, constrained bootstrap |
| $(S(\psi) - \hat{\mu}^S)/\hat{\sigma}^S$ | $S_{ls}^u$ | Location-and-scale adjusted $S(\psi)$, unconstrained bootstrap |
| $(S(\psi) - \tilde{\mu}^S)/\tilde{\sigma}^S$ | $S_{ls}^c$ | Location-and-scale adjusted $S(\psi)$, constrained bootstrap |
| $(R(\psi) - \hat{\mu}^R)/\hat{\sigma}^R$ | $R_{ls}^u$ | Location-and-scale adjusted $R(\psi)$, unconstrained bootstrap |
| $(R(\psi) - \tilde{\mu}^R)/\tilde{\sigma}^R$ | $R_{ls}^c$ | Location-and-scale adjusted $R(\psi)$, constrained bootstrap |

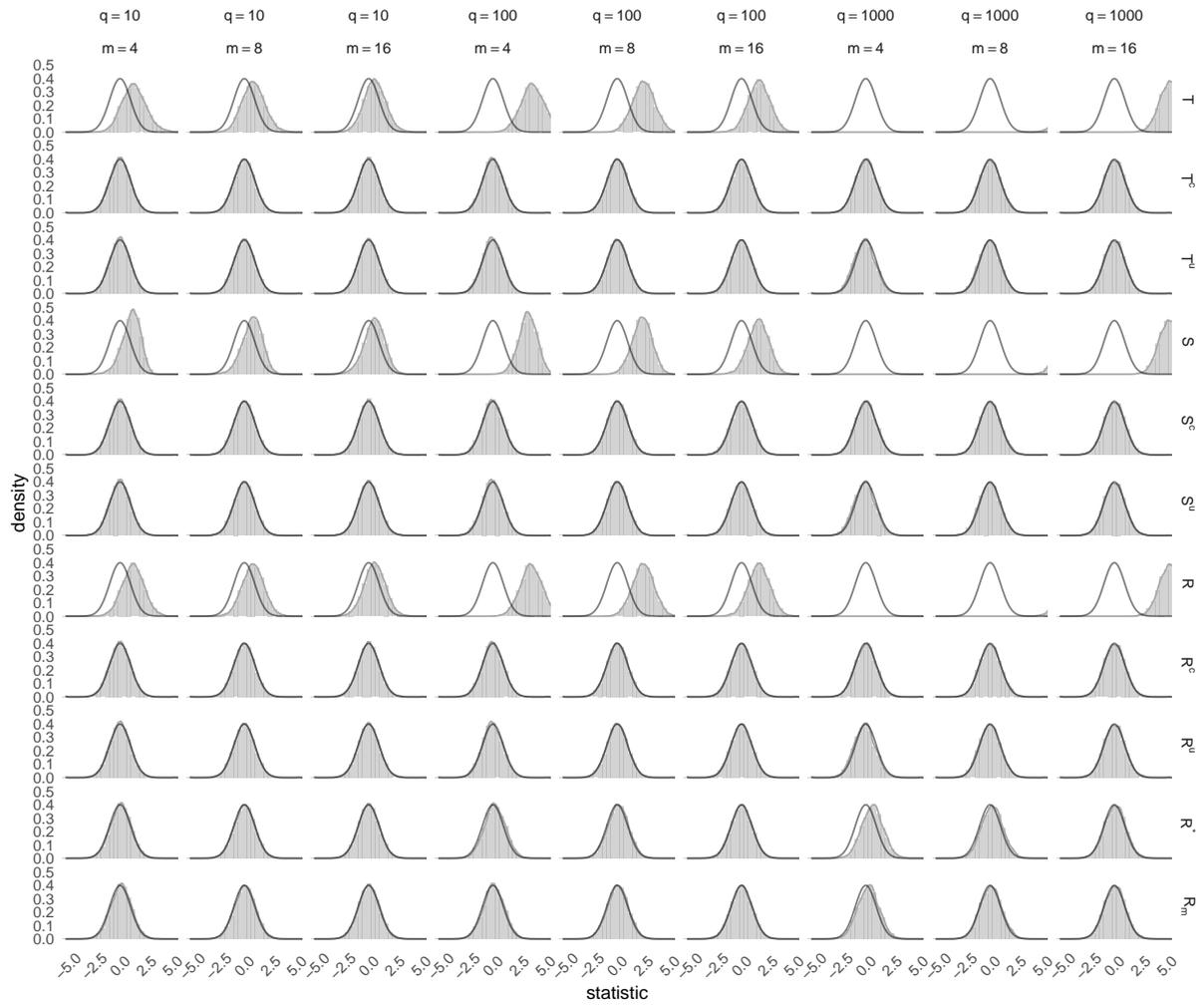

Figure 1: Gamma model. Estimated null distribution of statistics for the statistics in Table 1 for various combinations of $q$ and $m$. The N(0,1) density function is superimposed.

9

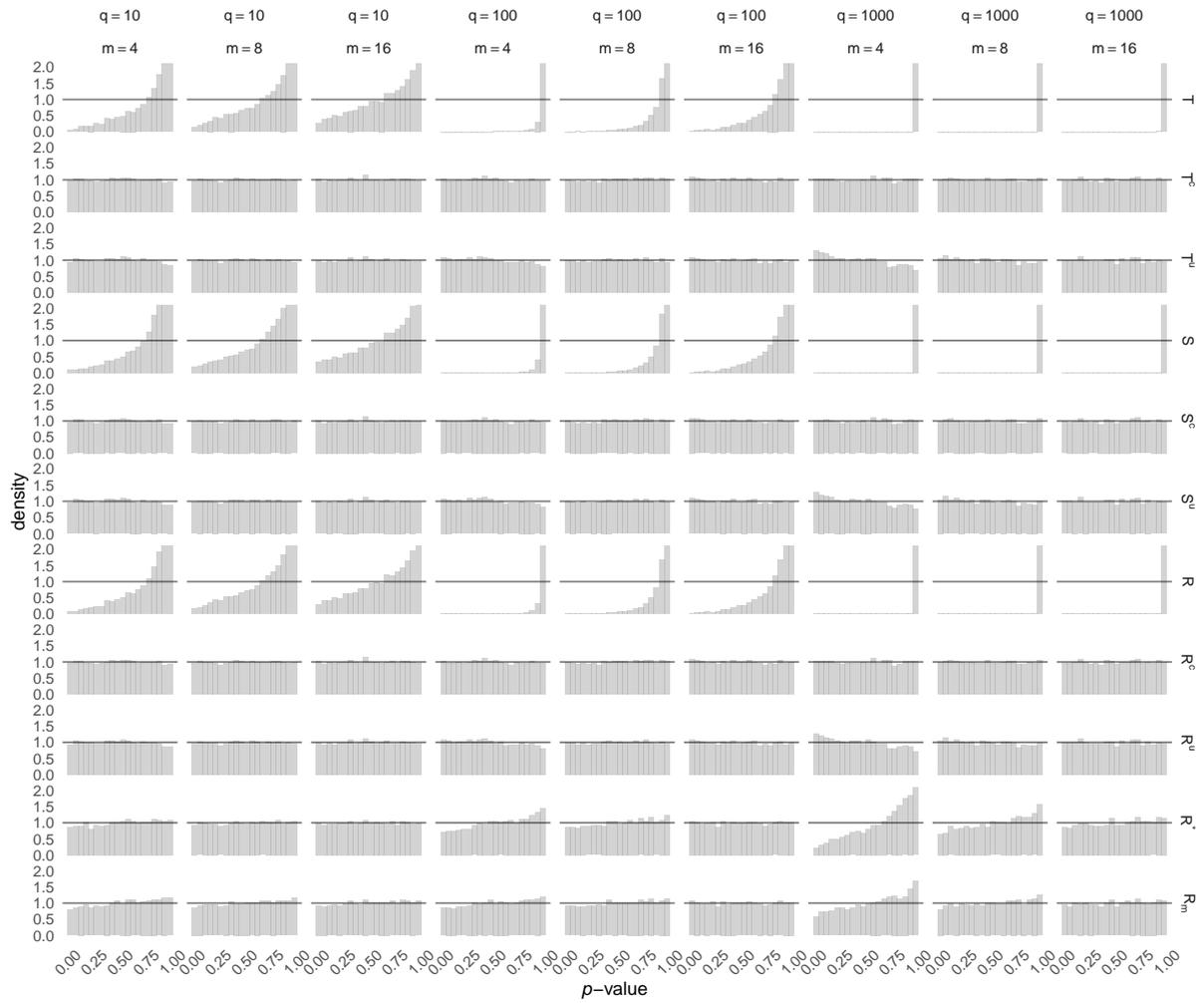

Figure 2: Gamma model. Estimated distribution of $p$-values for the statistics in Table 1 for various combinations of $q$ and $m$. The Uniform(0,1) density function is superimposed.

10

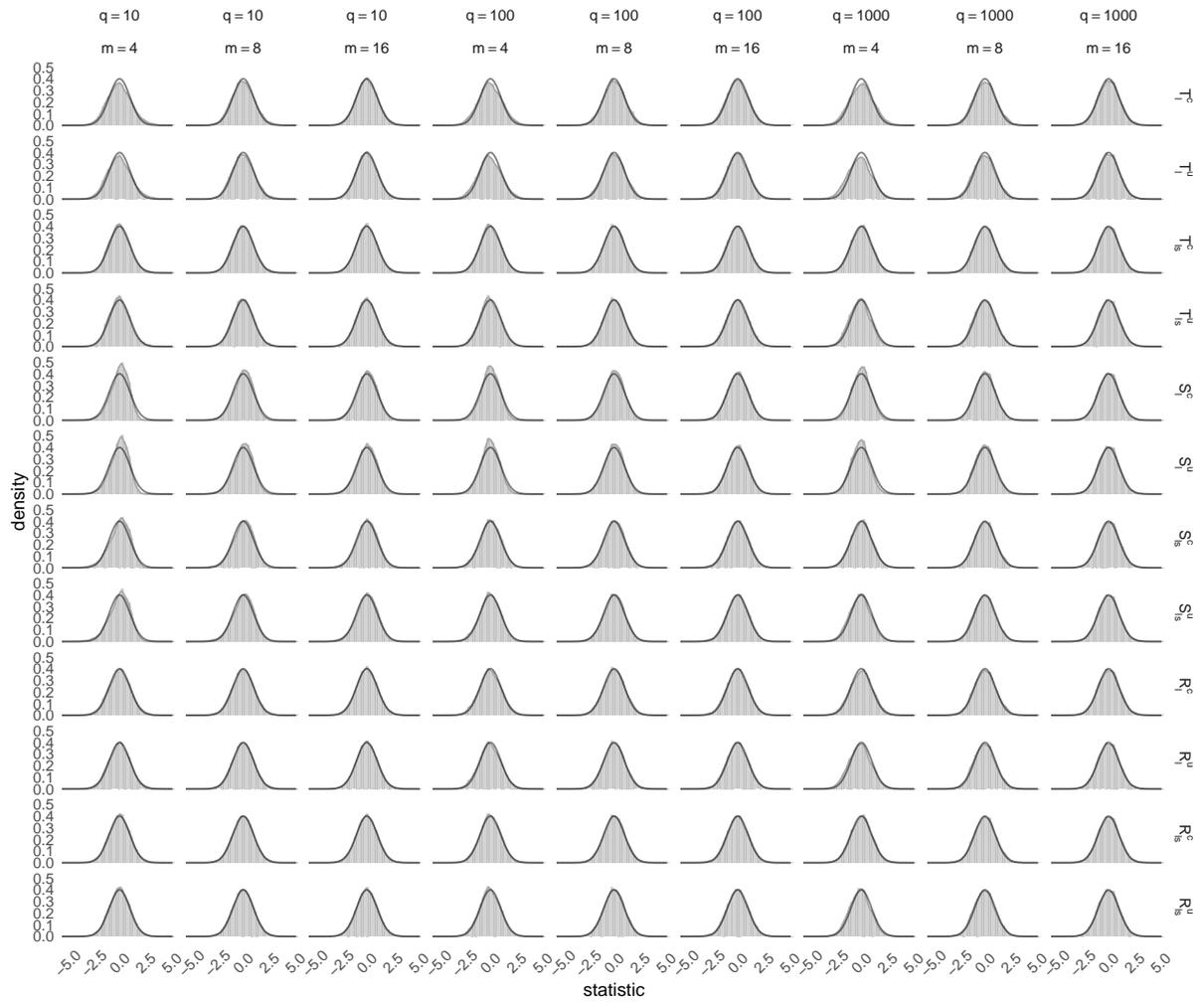

Figure 3: Gamma model. Estimated null distribution of statistics for the statistics in Table 2 for various combinations of $q$ and $m$. The N(0,1) density function is superimposed.

11

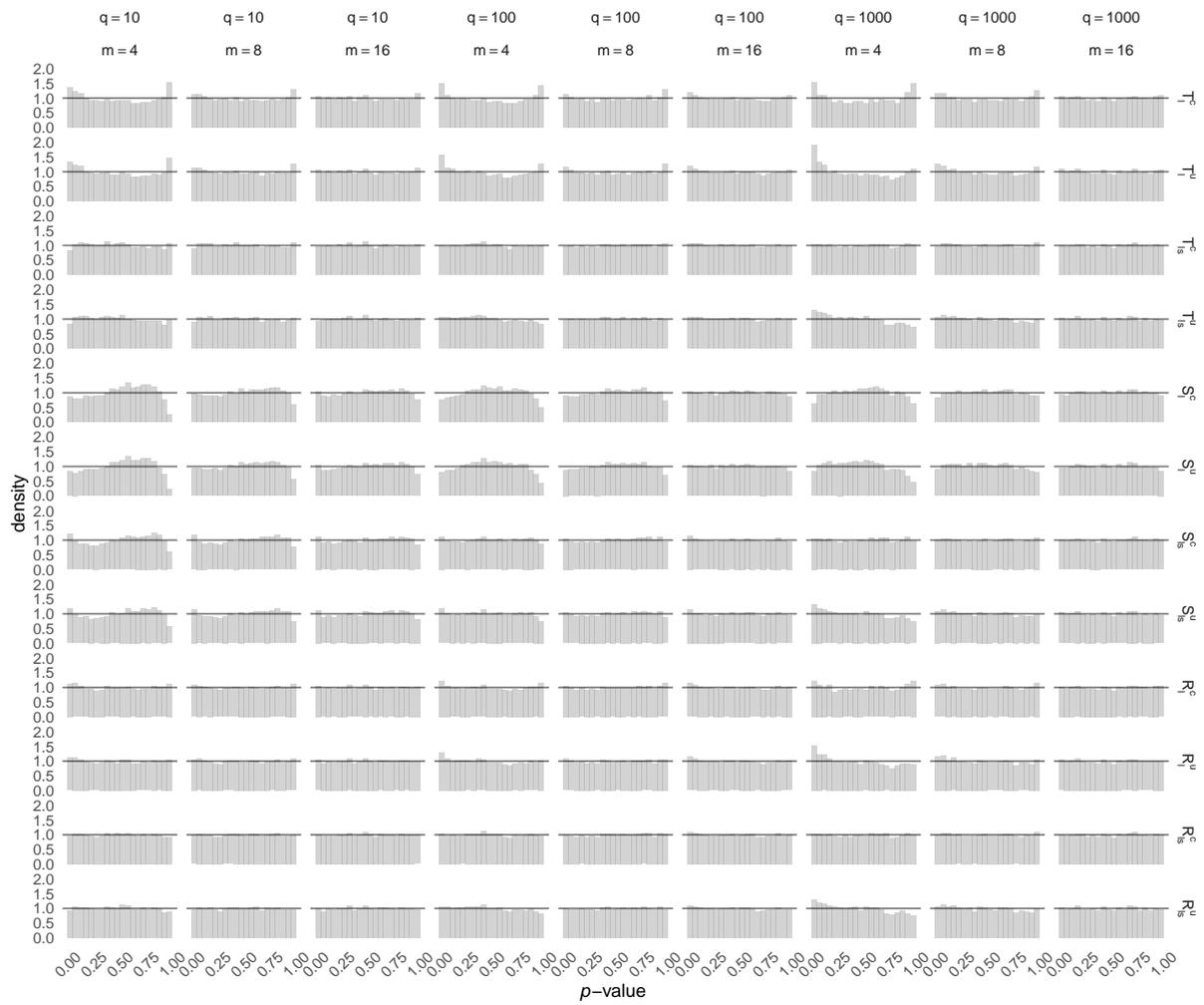

Figure 4: Gamma model. Estimated null distribution of $p$-values for the statistics in Table 2 for various combinations of $q$ and $m$. The Uniform(0,1) density function is superimposed.

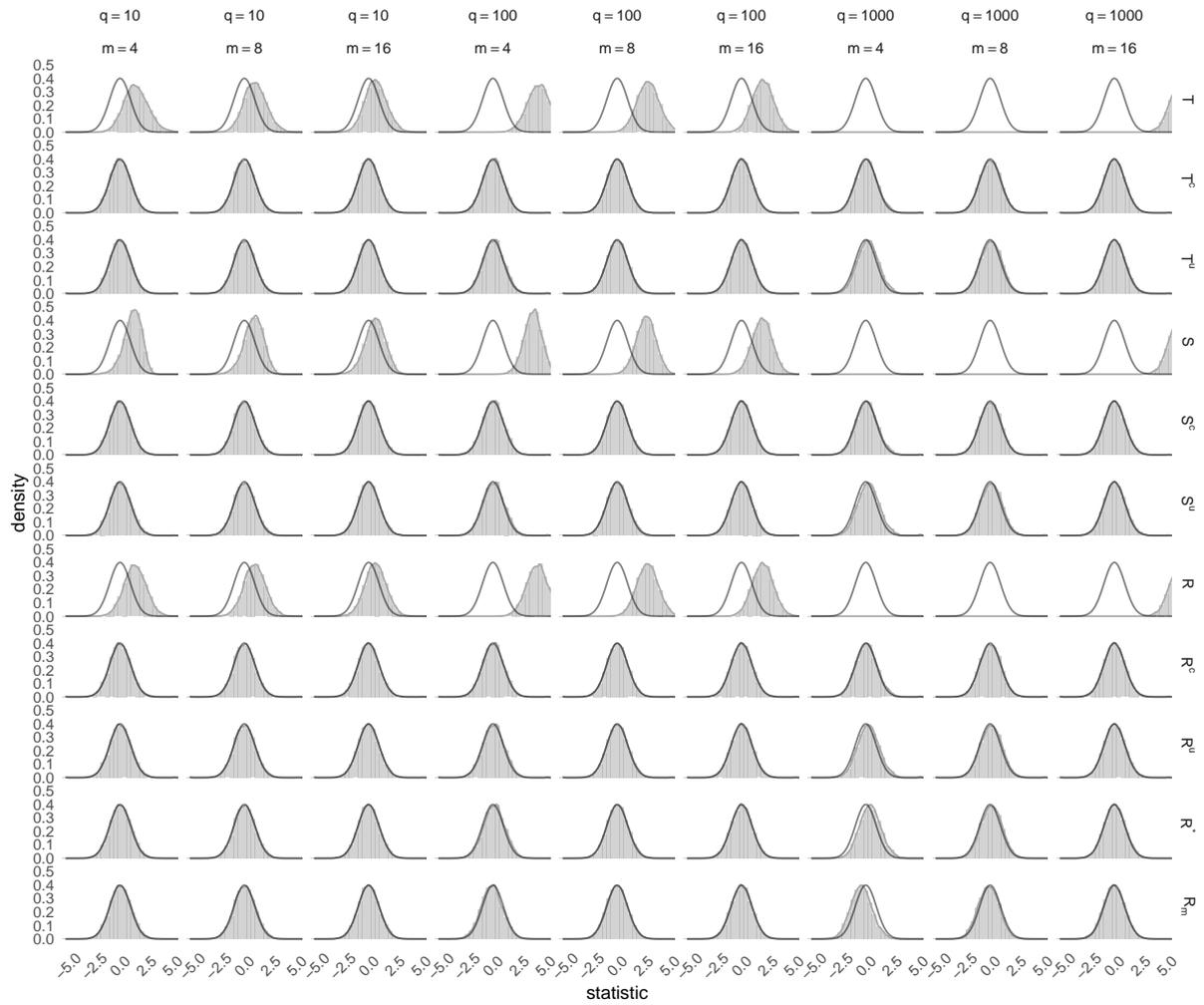

Figure 5: Beta model. Estimated null distribution of statistics for the statistics in Table 1 for various combinations of $q$ and $m$. The N(0,1) density function is superimposed.

13

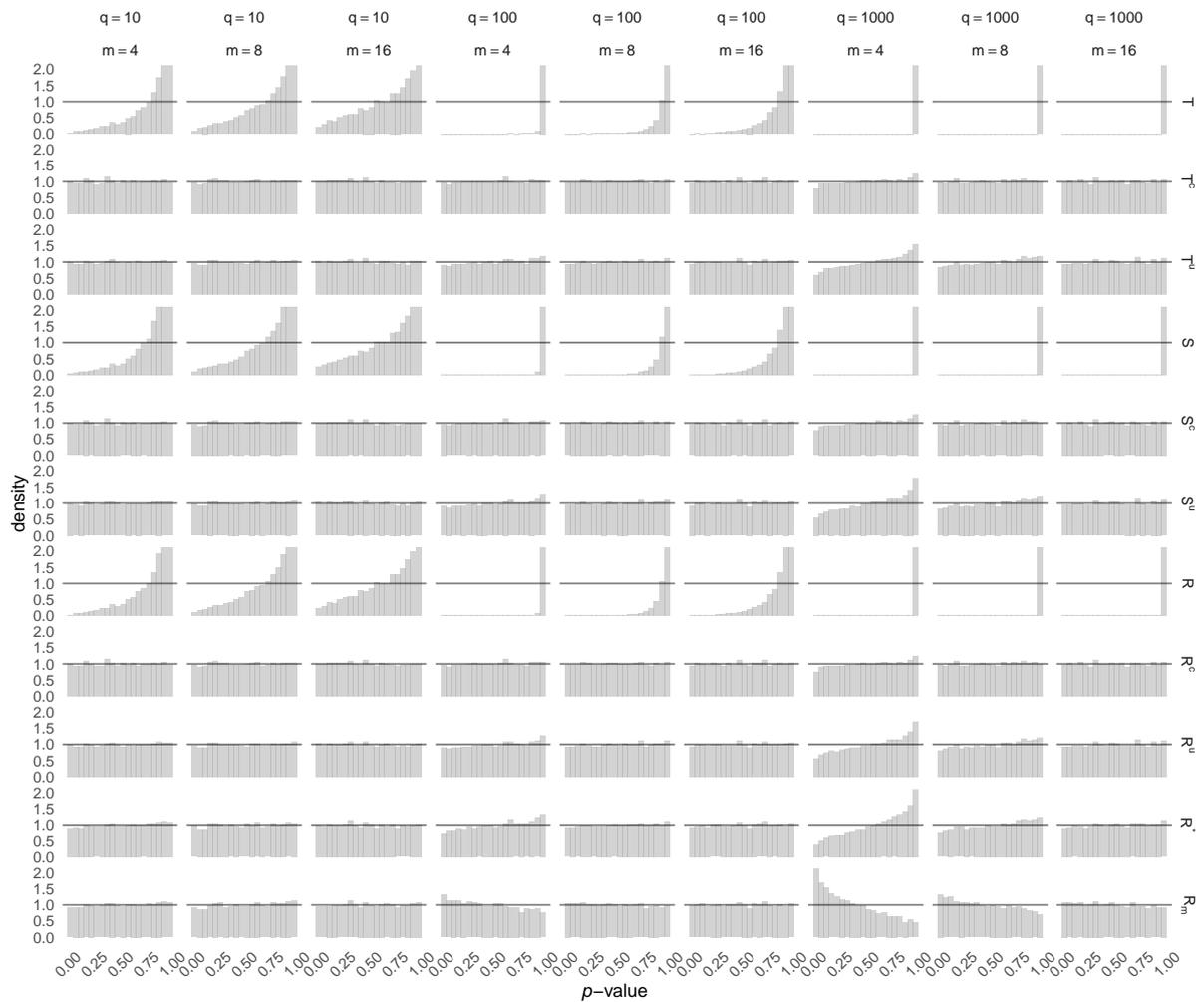

Figure 6: Beta model. Estimated distribution of $p$-values for the statistics in Table 1 for various combinations of $q$ and $m$. The Uniform(0,1) density function is superimposed.

14

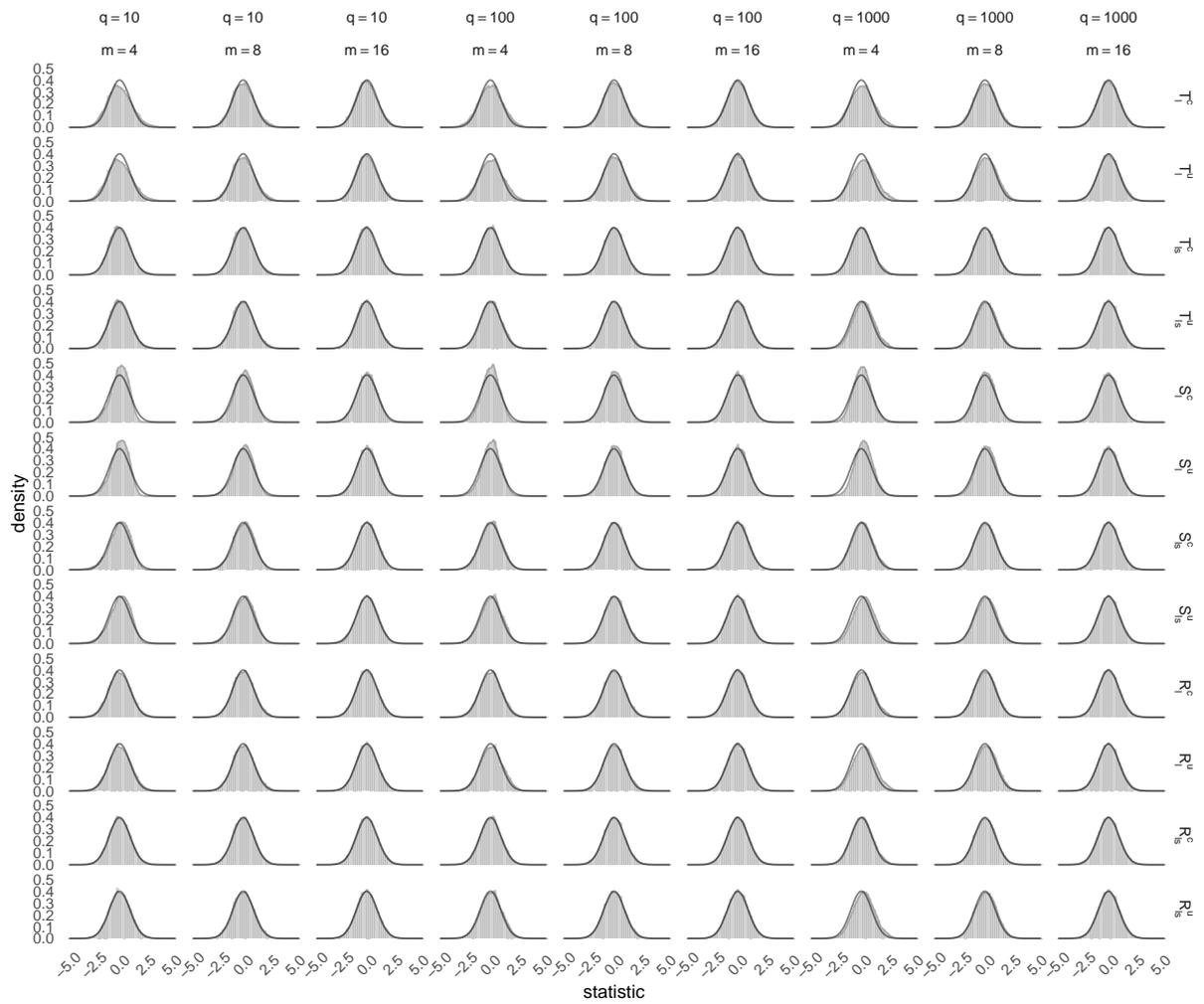

Figure 7: Beta model. Estimated null distribution of statistics for the statistics in Table 2 for various combinations of $q$ and $m$. The N(0,1) density function is superimposed.

15

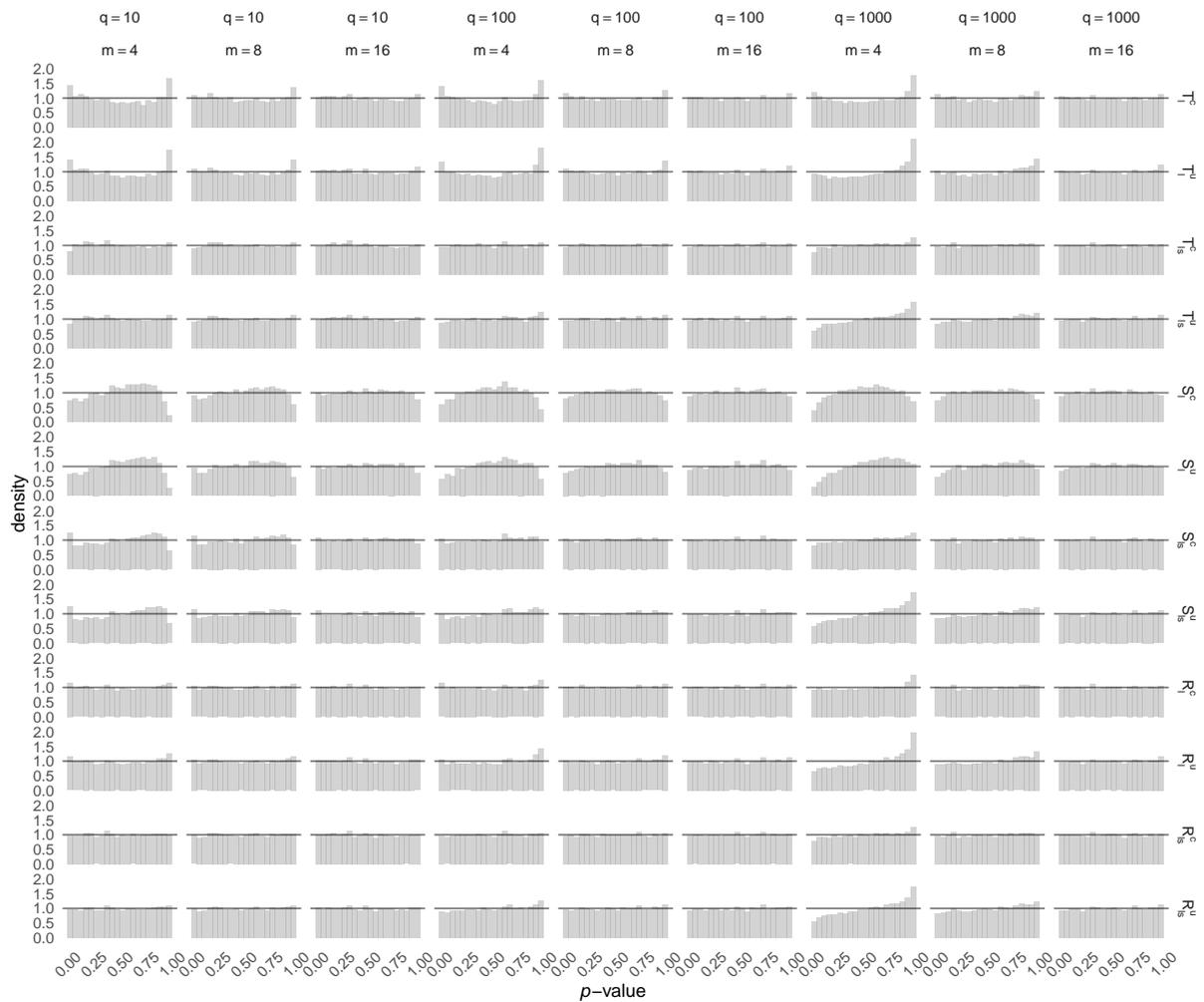

Figure 8: Beta model. Estimated null distribution of $p$-values for the statistics in Table 2 for various combinations of $q$ and $m$. The Uniform(0,1) density function is superimposed.

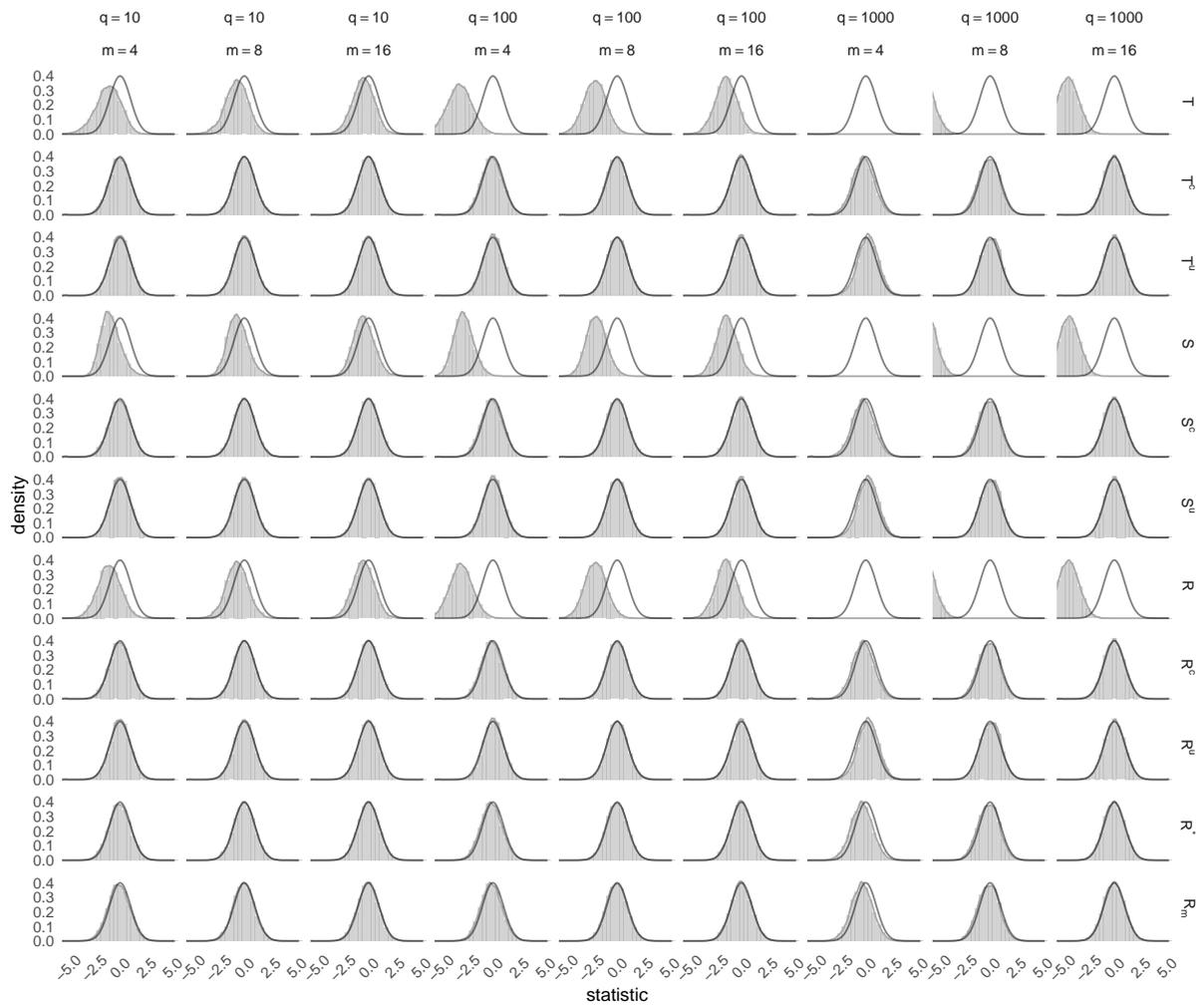

Figure 9: Curved exponential family. Estimated null distribution of statistics for the statistics in Table 1 for various combinations of $q$ and $m$. The N(0,1) density function is superimposed.

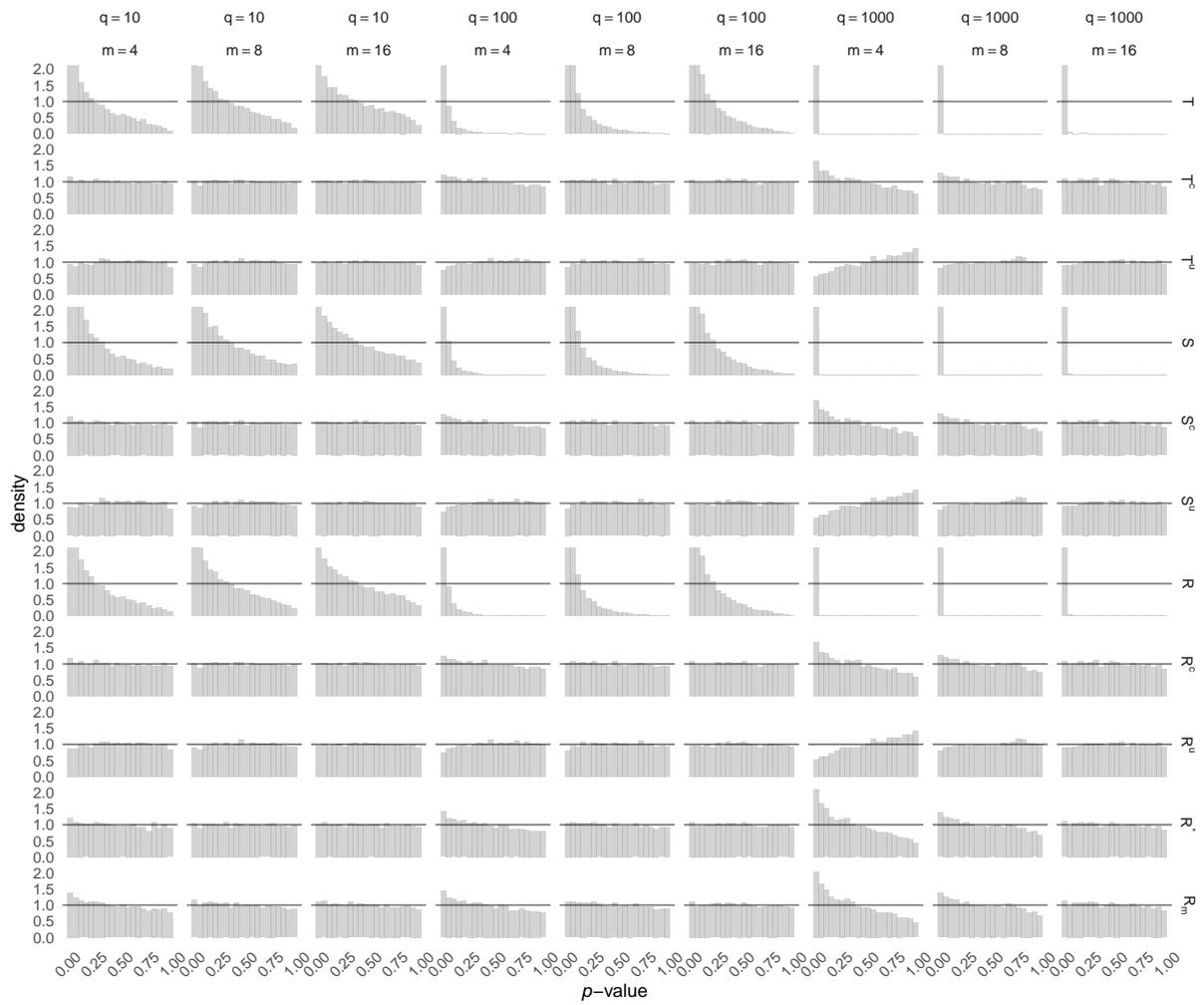

Figure 10: Curved exponential family. Estimated distribution of $p$-values for the statistics in Table 1 for various combinations of $q$ and $m$. The Uniform(0,1) density function is superimposed.

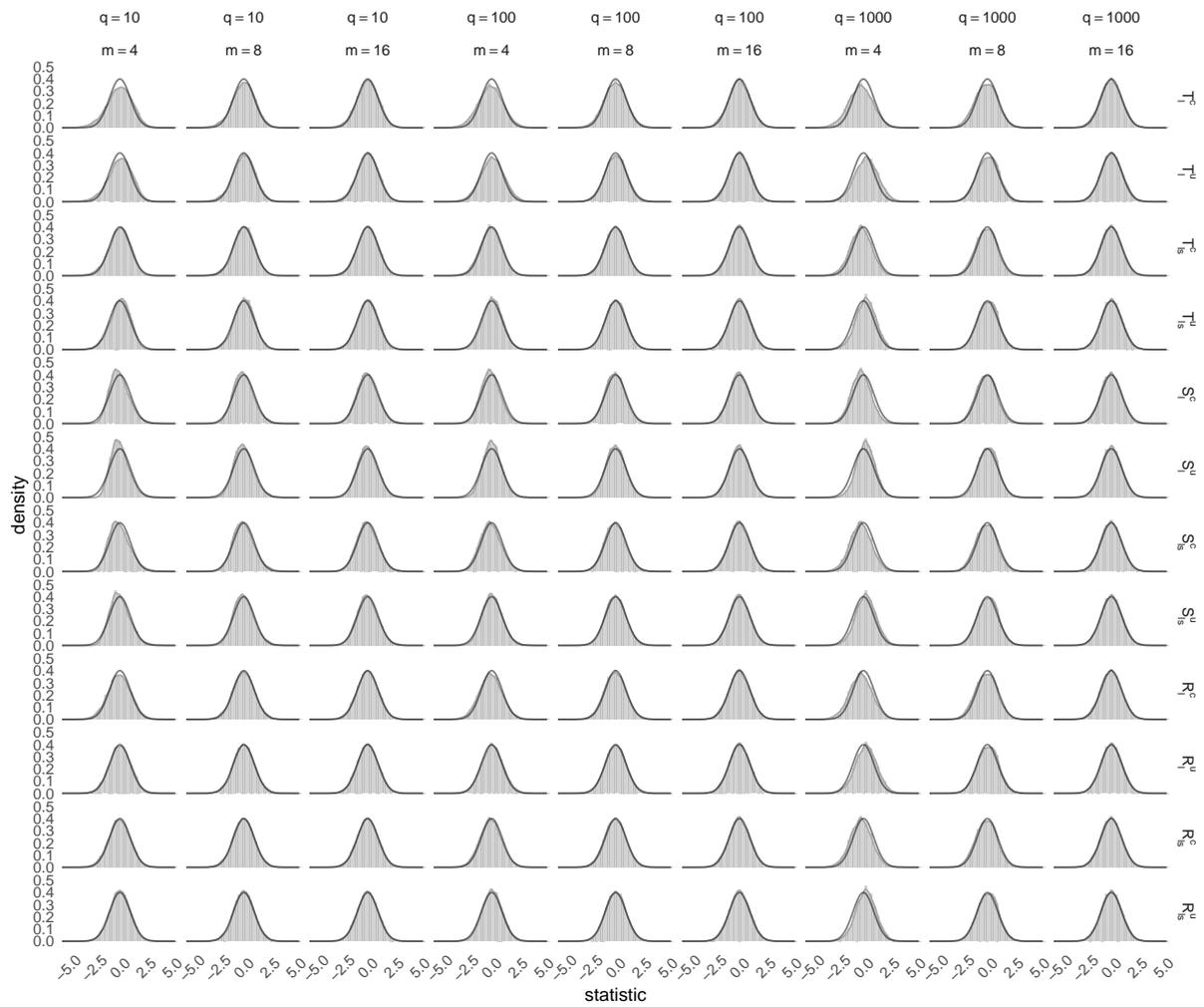

Figure 11: Curved exponential family. Estimated null distribution of statistics for the statistics in Table 2 for various combinations of $q$ and $m$. The N(0,1) density function is superimposed.

19

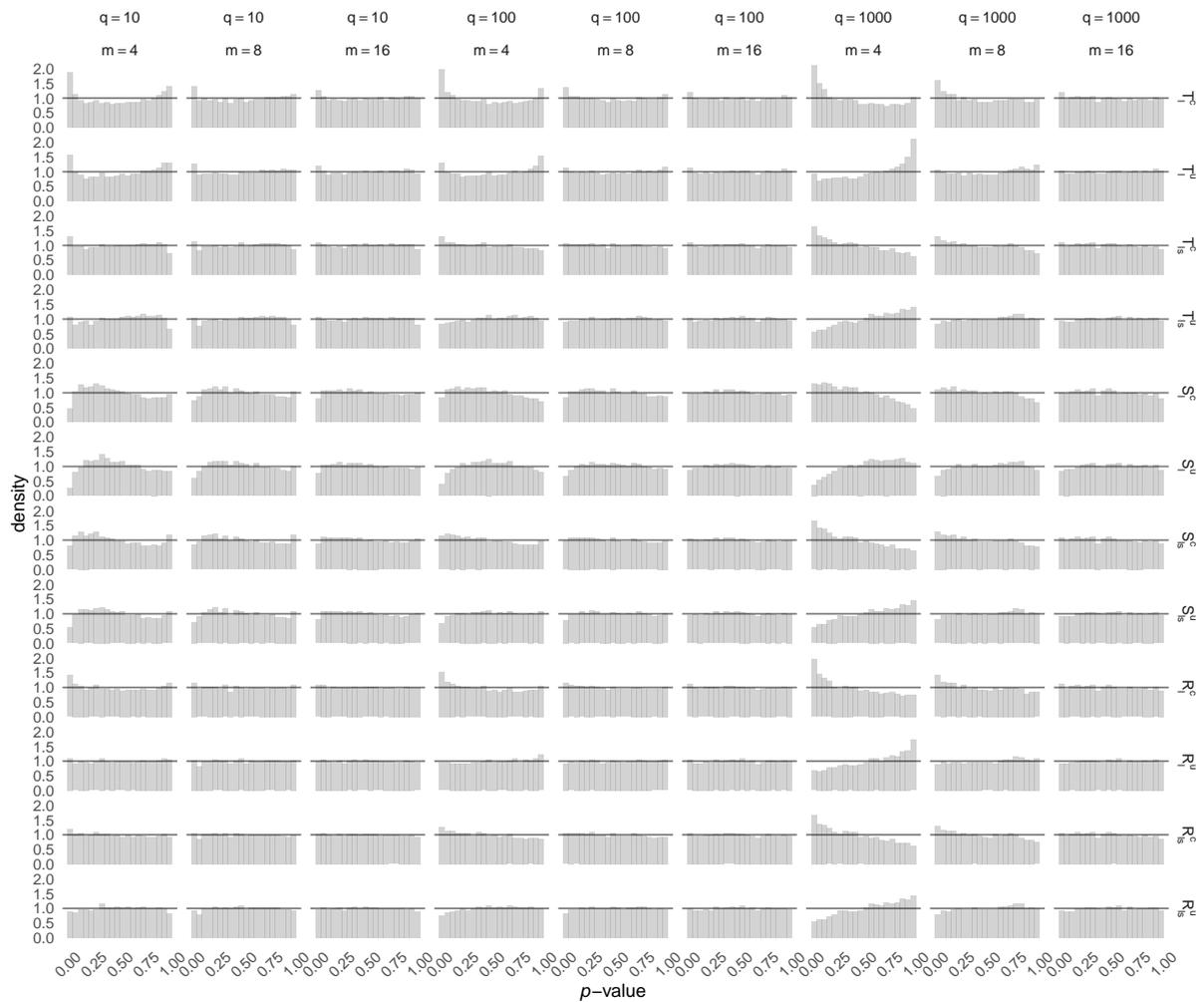

Figure 12: Curved exponential family. Estimated null distribution of $p$-values for the statistics in Table 2 for various combinations of $q$ and $m$. The Uniform(0,1) density function is superimposed.

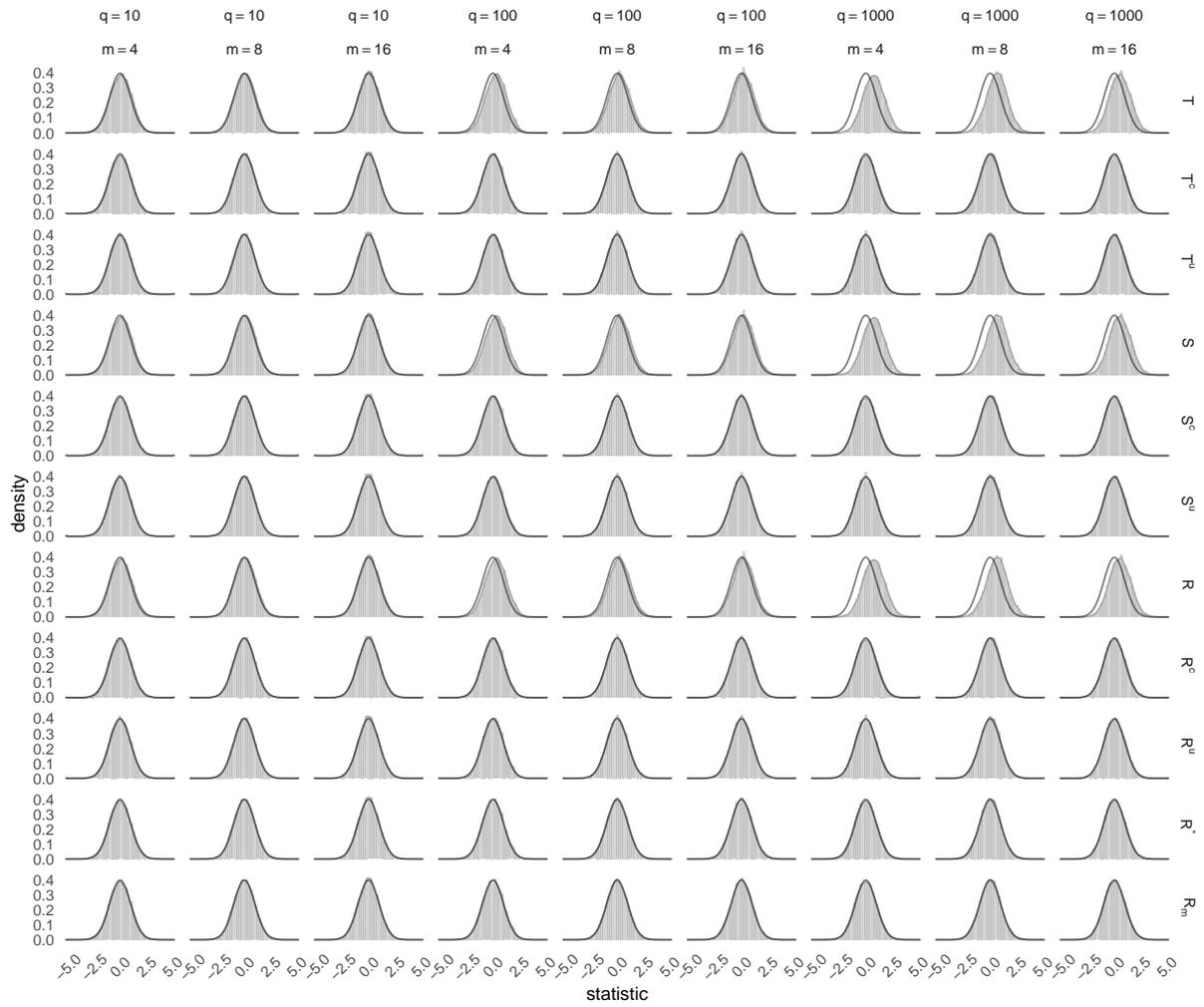

Figure 13: Truncated linear regression model. Estimated null distribution of statistics for the statistics in Table 1 for various combinations of $q$ and $m$. The N(0,1) density function is superimposed.

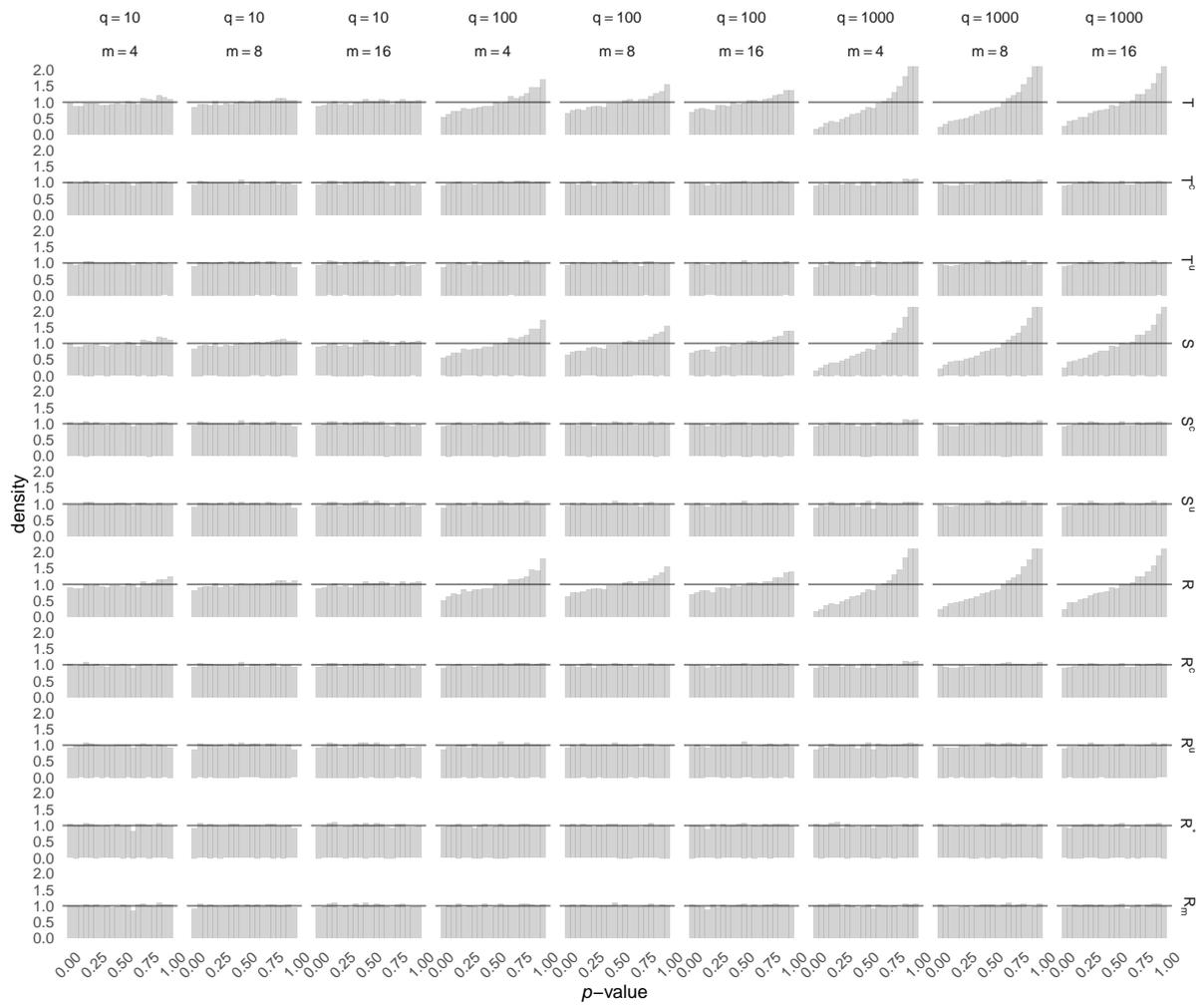

Figure 14: Truncated linear regression model. Estimated distribution of $p$-values for the statistics in Table 1 for various combinations of $q$ and $m$. The Uniform(0,1) density function is superimposed.

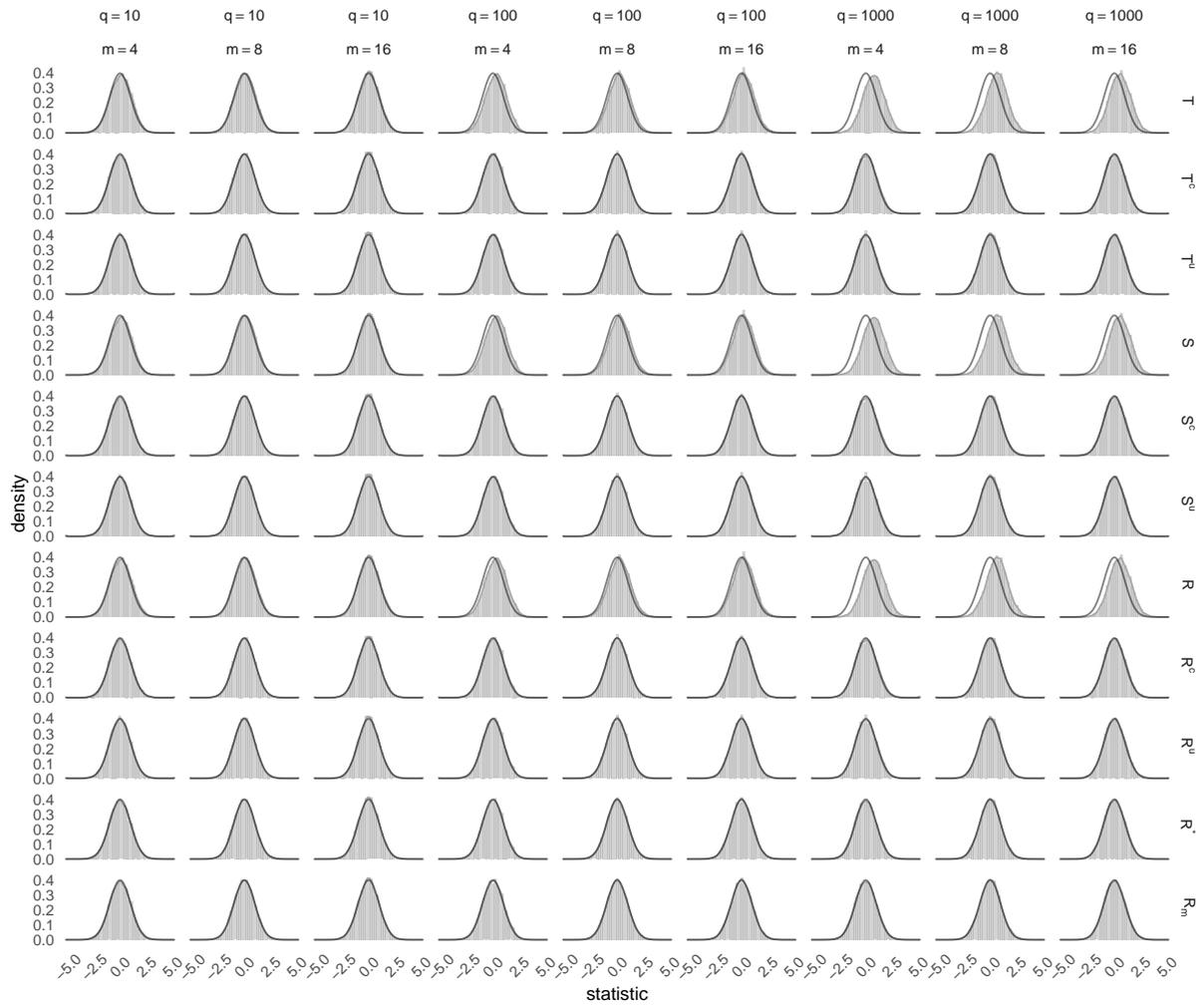

Figure 15: Truncated linear regression model. Estimated null distribution of statistics for the statistics in Table 2 for various combinations of $q$ and $m$. The N(0,1) density function is superimposed.

23

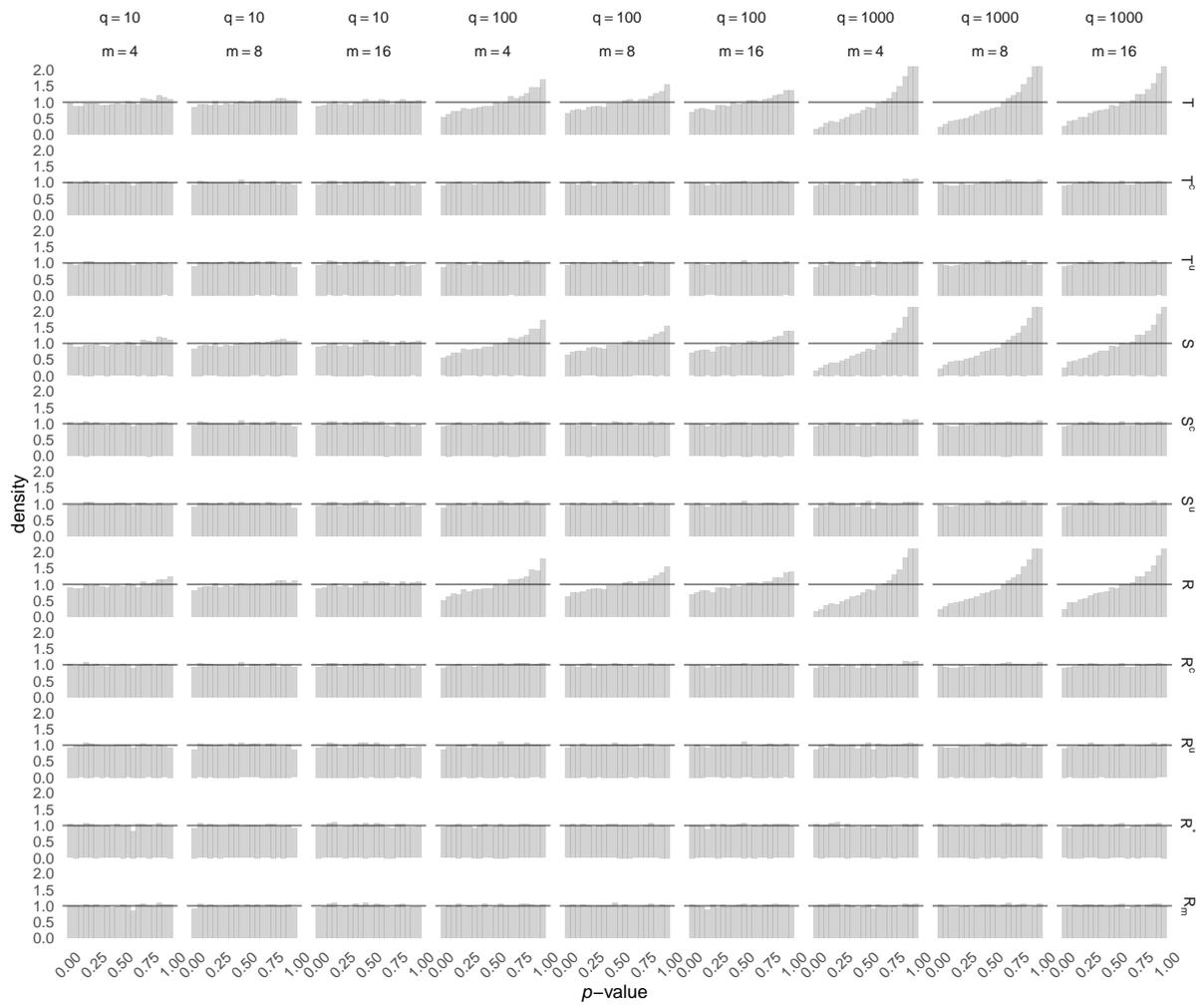

Figure 16: Truncated linear regression model. Estimated null distribution of $p$-values for the statistics in Table 2 for various combinations of $q$ and $m$. The Uniform(0,1) density function is superimposed.

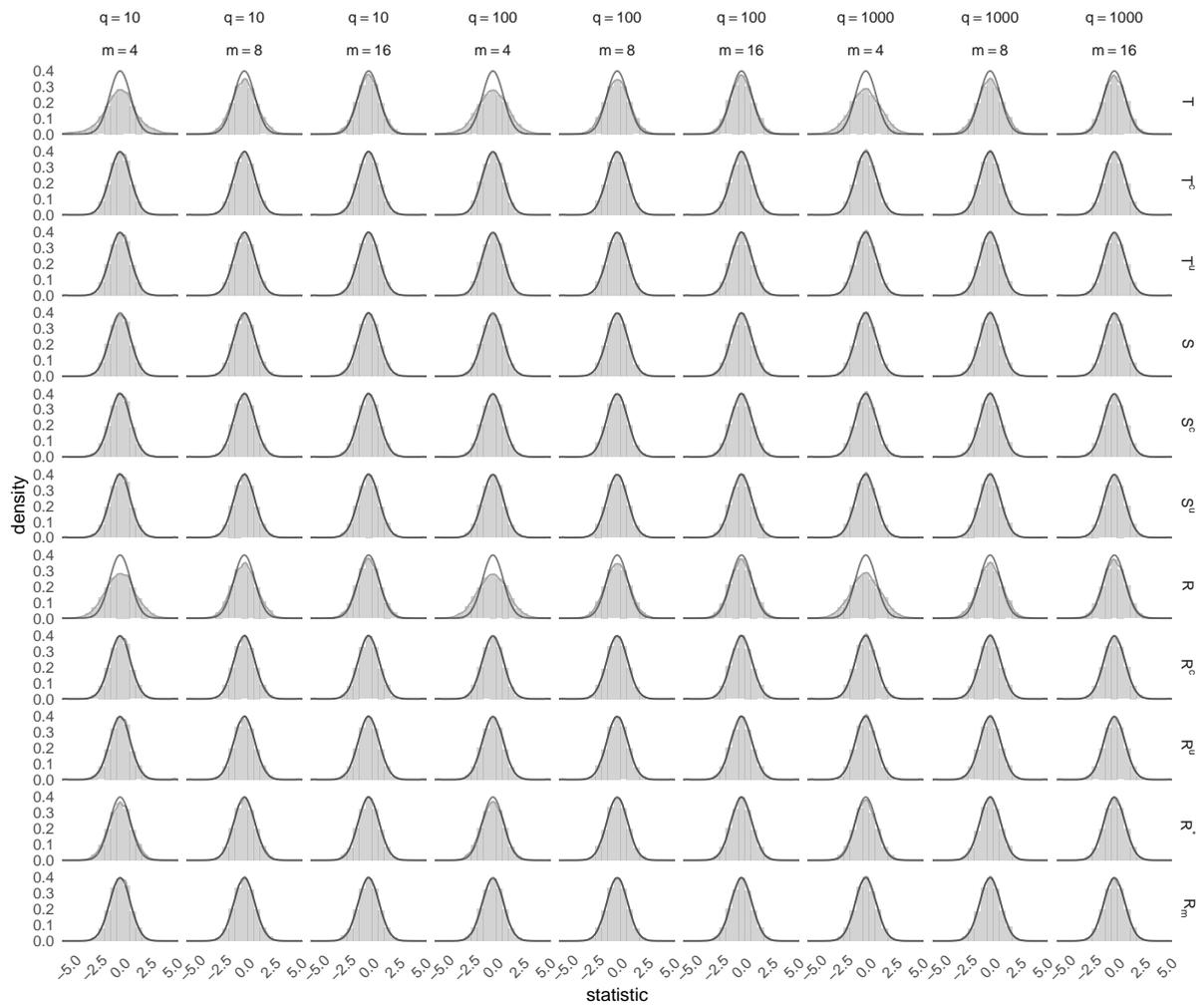

Figure 17: Behrens-Fisher model. Estimated null distribution of statistics for the statistics in Table 1 for various combinations of $q$ and $m$. The N(0,1) density function is superimposed.

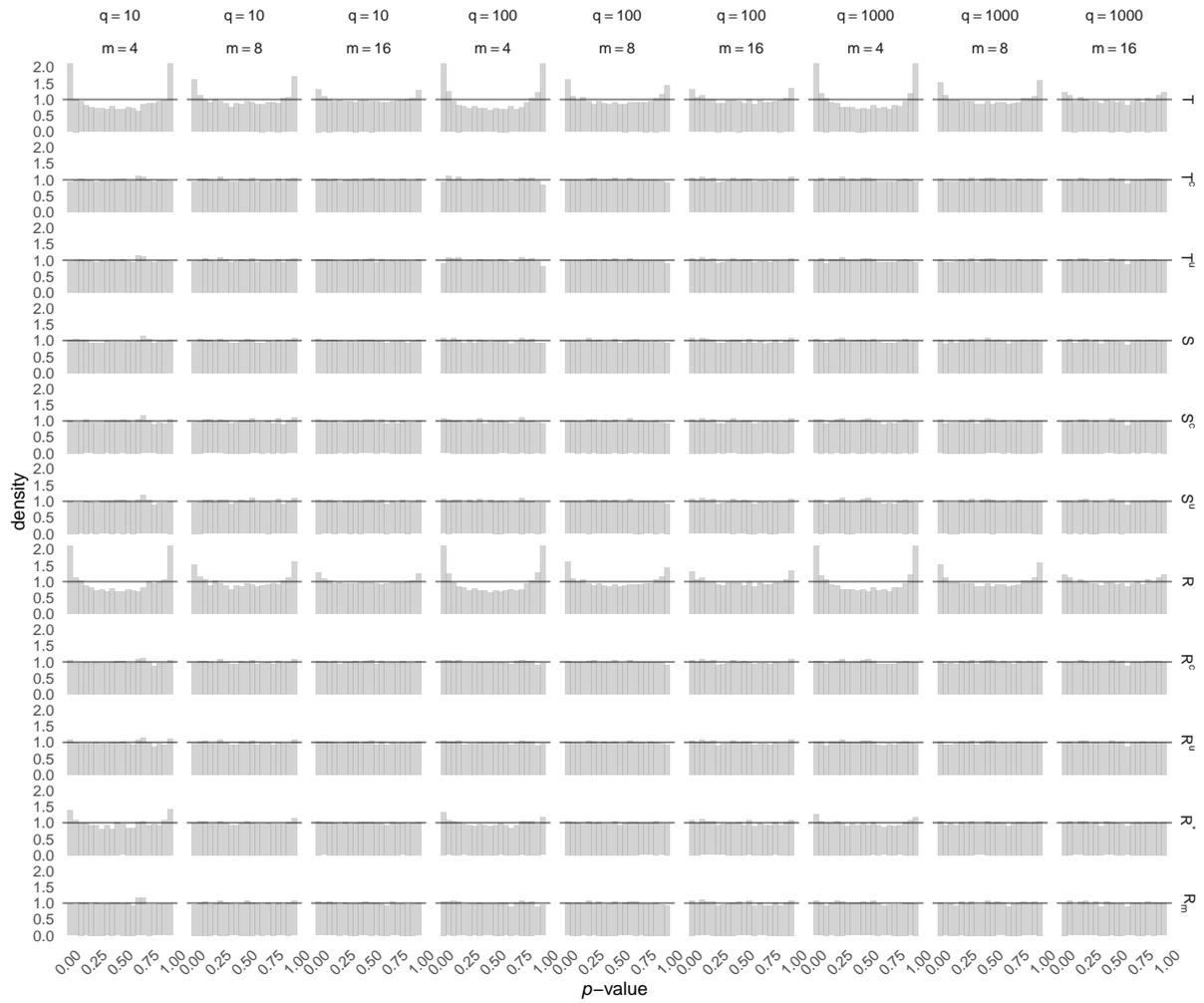

Figure 18: Behrens-Fisher model. Estimated distribution of $p$-values for the statistics in Table 1 for various combinations of $q$ and $m$. The Uniform(0,1) density function is superimposed.

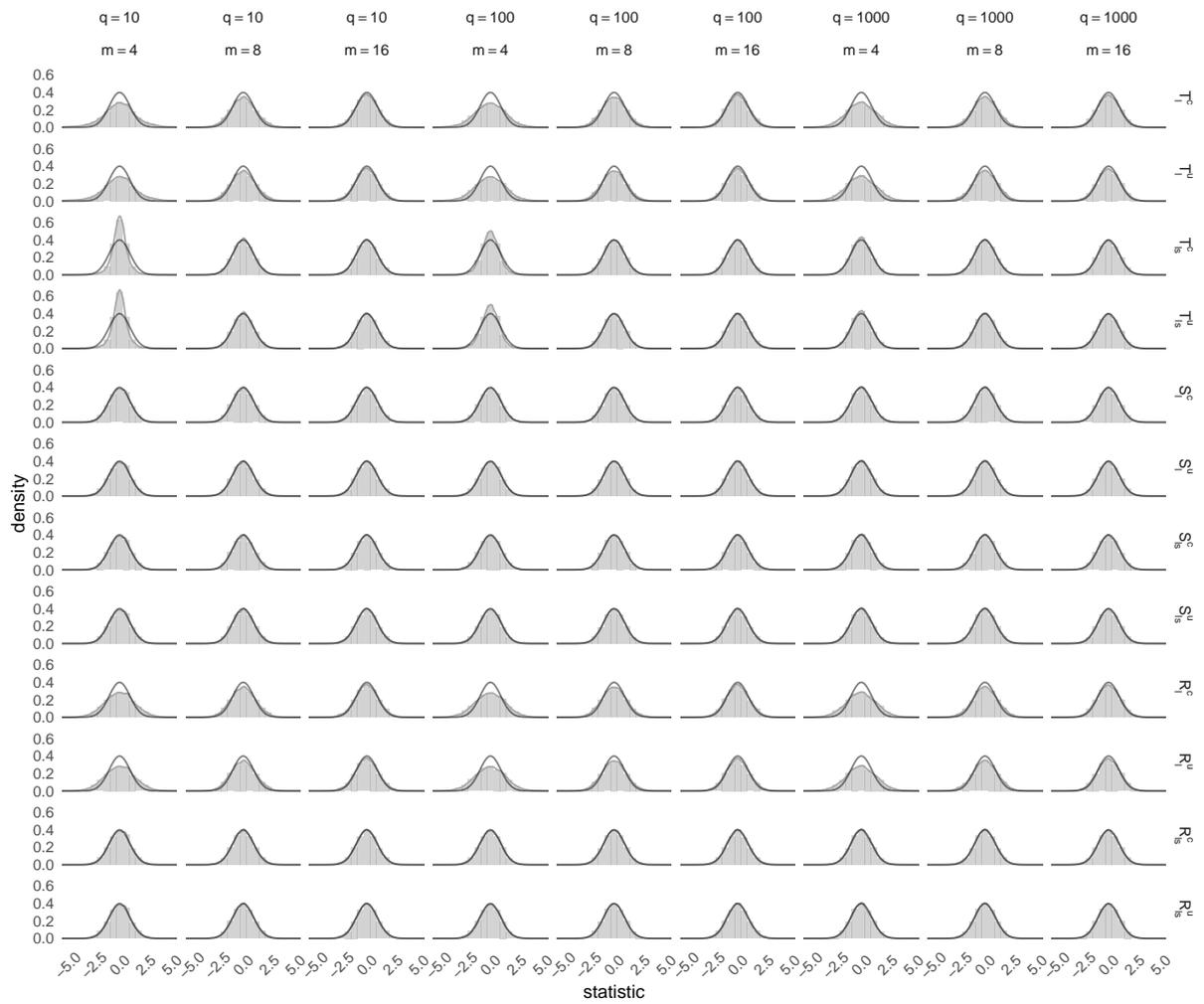

Figure 19: Behrens-Fisher model. Estimated null distribution of statistics for the statistics in Table 2 for various combinations of $q$ and $m$. The N(0,1) density function is superimposed.

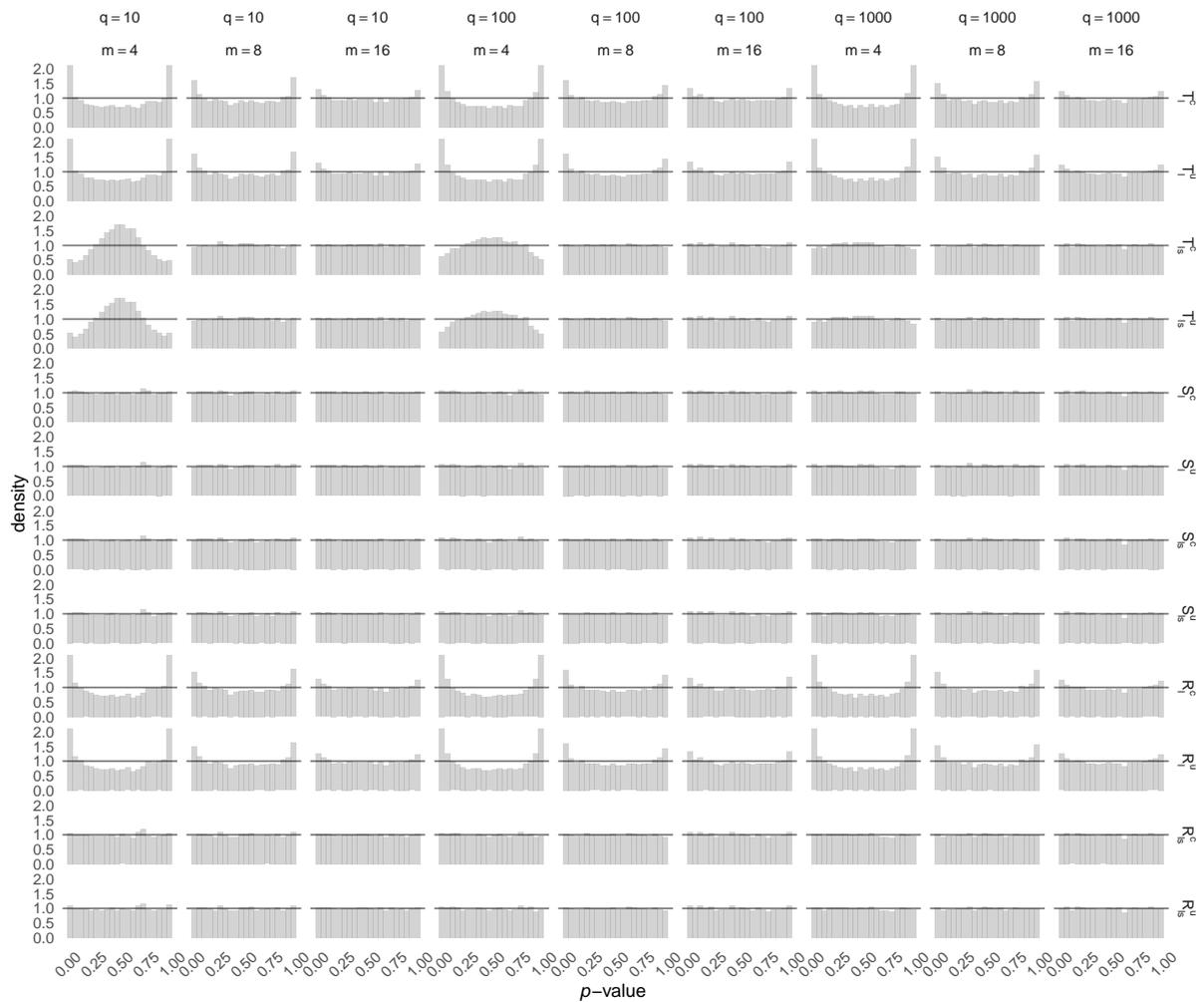

Figure 20: Behrens-Fisher model. Estimated null distribution of $p$-values for the statistics in Table 2 for various combinations of $q$ and $m$. The Uniform(0,1) density function is superimposed.

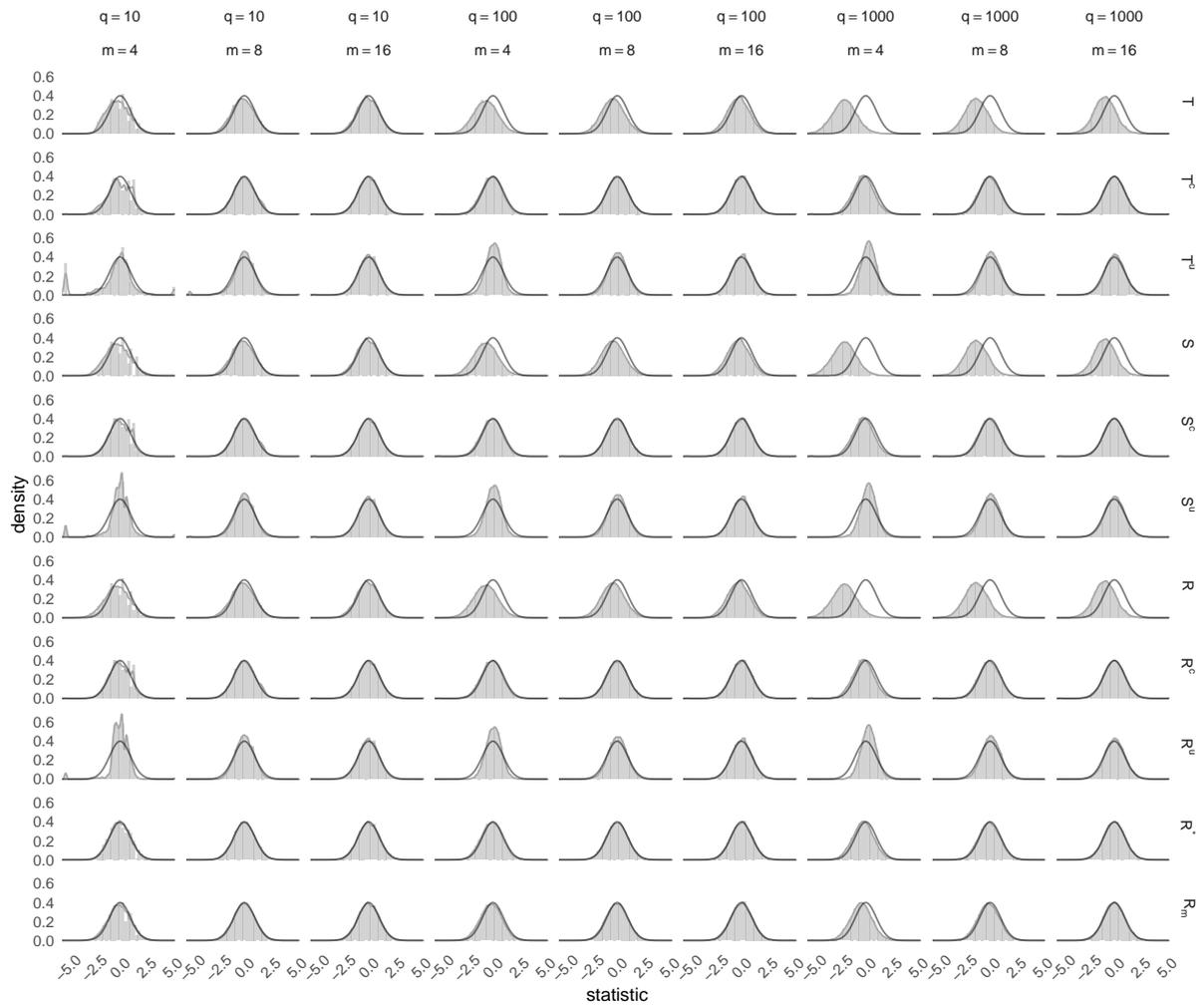

Figure 21: Logistic regression model. Estimated null distribution of statistics for the statistics in Table 1 for various combinations of $q$ and $m$. The N(0,1) density function is superimposed.

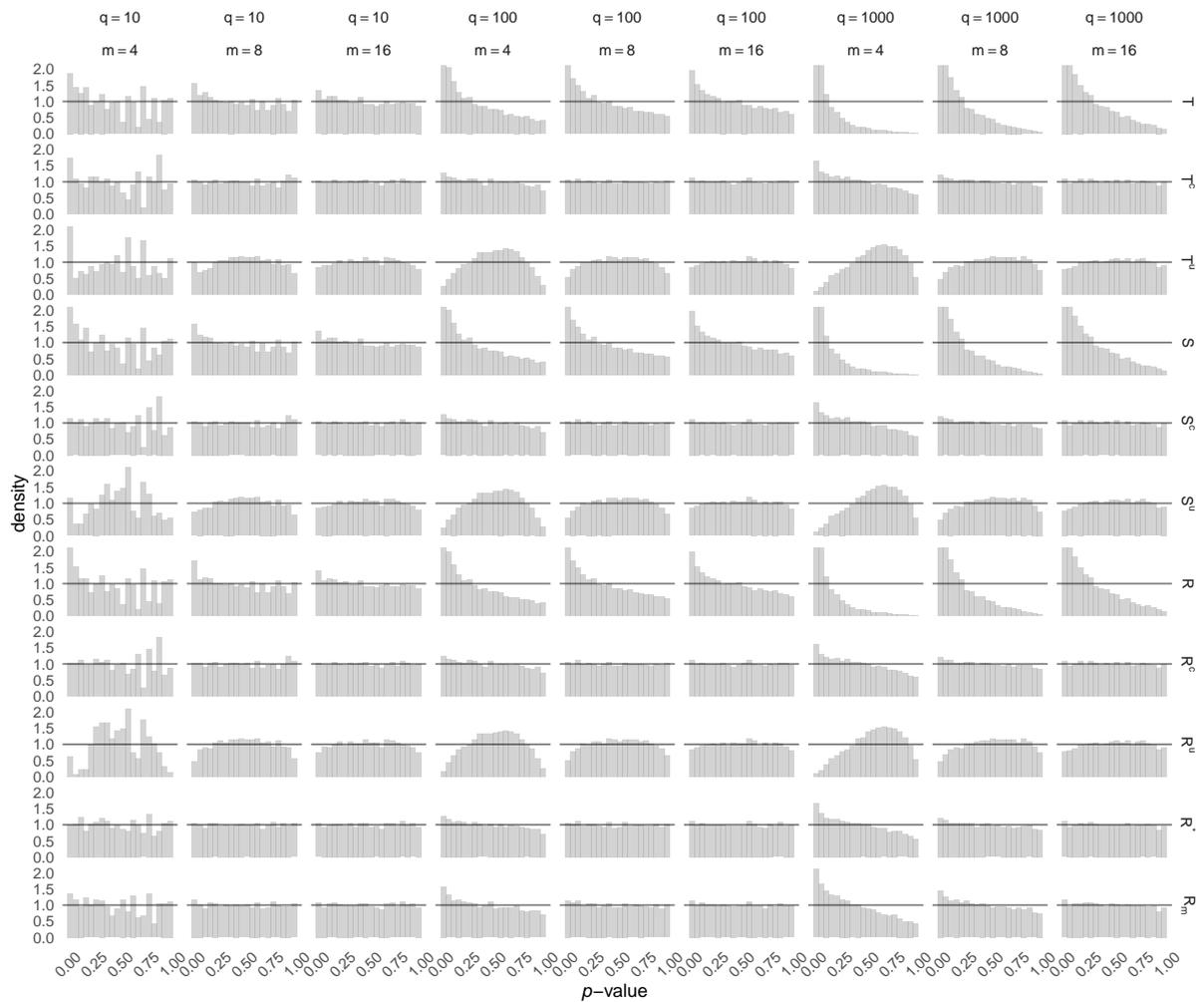

Figure 22: Logistic regression model. Estimated distribution of $p$-values for the statistics in Table 1 for various combinations of $q$ and $m$. The Uniform(0,1) density function is superimposed.

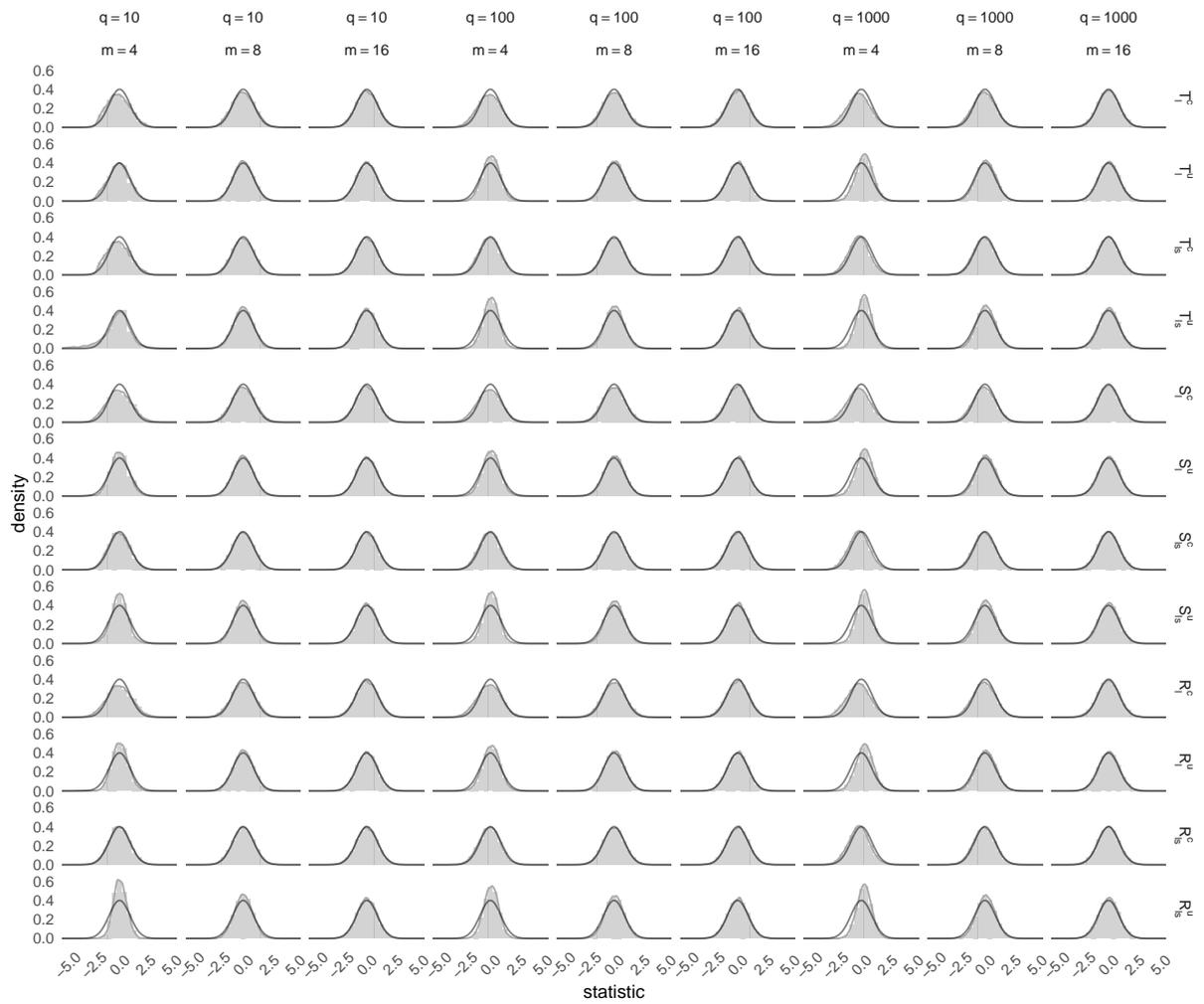

Figure 23: Logistic regression model. Estimated null distribution of statistics for the statistics in Table 2 for various combinations of $q$ and $m$. The N(0,1) density function is superimposed.

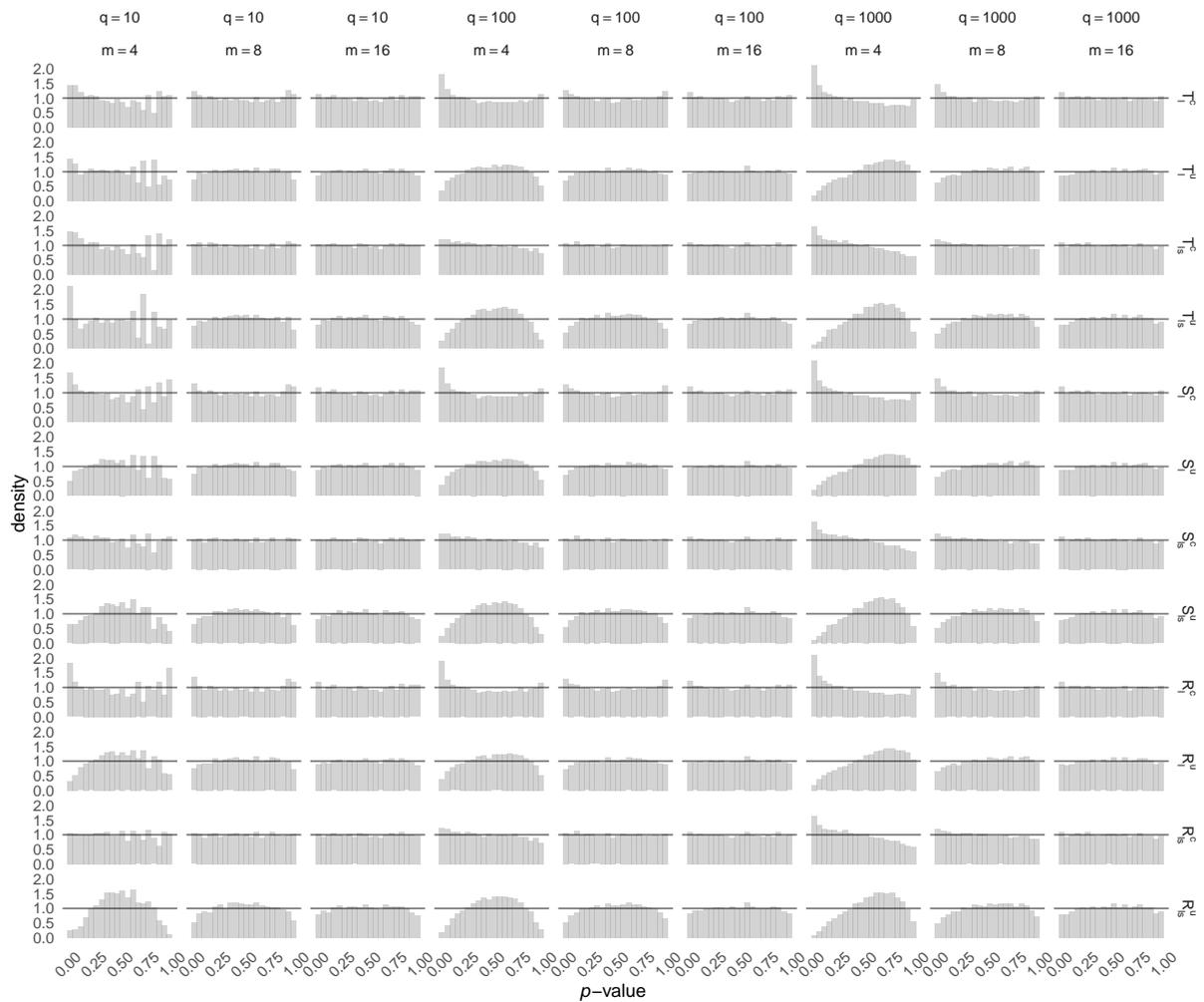

Figure 24: Logistic regression model. Estimated null distribution of $p$-values for the statistics in Table 2 for various combinations of $q$ and $m$. The Uniform(0,1) density function is superimposed.

32

# S3 Empirical tail probabilities

Tables 3-11 report empirical tail probabilities for the statistics in Table 1 of the main text, for all models and all simulation settings. The rows of each table are subdivided into 6 different blocks, each one corresponding to a different model.



\end{document}